\documentclass[11pt]{amsart}
\usepackage{amsmath}
\usepackage{amsxtra}
\usepackage{amscd}
\usepackage{amsthm}
\usepackage{amsfonts}
\usepackage{amssymb}
\usepackage{eucal}

\newtheorem{cor}[subsubsection]{Corollary}
\newtheorem{lem}[subsubsection]{Lemma}
\newtheorem{prop}[subsubsection]{Proposition}

\newtheorem{thm}[subsubsection]{Theorem}

\theoremstyle{definition}
\theoremstyle{remark}

\newcommand{\nc}{\newcommand}
\nc{\renc}{\renewcommand}
\nc{\ssec}{\subsection}
\nc{\sssec}{\subsubsection}
\nc{\on}{\operatorname}

\nc\ol{\overline}
\nc\wt{\widetilde}
\nc\tboxtimes{\wt{\boxtimes}}
\nc{\alp}{\alpha}

\nc{\ZZ}{{\mathbb Z}}
\nc{\NN}{{\mathbb N}}
\nc{\CC}{{\mathbb C}}
\nc{\OO}{{\mathbb O}}
\renc{\SS}{{\mathbb S}}
\nc{\DD}{{\mathbb D}}
\nc{\GG}{{\mathbb G}}

\nc{\Fq}{{\mathbb F}_q}
\nc{\Fqb}{\ol{{\mathbb F}_q}}
\nc{\Ql}{\ol{{\mathbb Q}_\ell}}
\nc{\id}{\text{id}}
\nc\X{\mathcal X}

\nc{\Hom}{\on{Hom}}
\nc{\Lie}{\on{Lie}}
\nc{\Loc}{\on{Loc}}
\nc{\Pic}{\on{Pic}}
\nc{\Bun}{\on{Bun}}
\nc{\IC}{\on{IC}}
\nc{\Aut}{\on{Aut}}
\nc{\rk}{\on{rk}}
\nc{\Sh}{\on{Sh}}
\nc{\IrrelSh}{\on{IrrelSh}}
\nc{\Perv}{\on{Perv}}
\nc{\pos}{{\on{pos}}}
\nc{\Conv}{\on{Conv}}
\nc{\Sph}{\on{Sph}}
\nc{\Sym}{\on{Sym}}
\nc{\BunBb}{\overline{\Bun}_B}
\nc{\Buno}{\overset{o}{\Bun}}
\nc{\BunPb}{{\overline{\Bun}_P}}
\nc{\BunBM}{\overline{\Bun}_{B(M)}}
\nc{\BunPbw}{{\widetilde{\Bun}_P}}
\nc{\BunBP}{\widetilde{\Bun}_{B,P}}
\nc{\GUb}{\overline{G/U}}
\nc{\GUPb}{\overline{G/U(P)}}

\nc{\Hhom}{\underline{\on{Hom}}}
\nc\syminfty{\on{Sym}^{\infty}}
\nc\lal{\ol{\lambda}}
\nc\xl{\ol{x}}
\nc\thl{\ol{\theta}}
\nc\nul{\ol{\nu}}
\nc\mul{\ol{\mu}}
\nc\Sum\Sigma
\nc{\oX}{\overset{o}{X}{}}
\nc{\M}{{\mathcal M}}
\nc{\N}{{\mathcal N}}
\nc{\F}{{\mathcal F}}
\nc{\D}{{\mathcal D}}
\nc{\Q}{{\mathcal Q}}
\nc{\Y}{{\mathcal Y}}
\nc{\G}{{\mathcal G}}
\nc{\E}{{\mathcal E}}
\nc{\CalC}{{\mathcal C}}
\nc\Dh{\widehat{\D}}

\nc{\C}{{\mathcal C}}
\nc{\K}{{\mathcal K}}
\renewcommand{\H}{{\mathcal H}}

\nc{\T}{{\mathcal T}}
\nc{\V}{{\mathcal V}}
\renc{\P}{{\mathcal P}}
\nc{\A}{{\mathcal A}}
\nc{\B}{{\mathcal B}}
\nc{\U}{{\mathcal U}}

\nc{\Gr}{\on{Gr}}

\nc{\frn}{{\check{\mathfrak u}(P)}}
\nc{\p}{\mathfrak p}
\nc{\q}{\mathfrak q}
\nc\f{{\mathfrak f}}

\nc{\qo}{{\mathfrak q}}
\nc{\po}{{\mathfrak p}}
\nc{\s}{{\mathfrak s}}
\nc\w{\text{w}}

\nc\mathi\iota
\nc\Spec{\on{Spec}}
\nc\Mod{\on{Mod}}
\nc{\tw}{\widetilde{\mathfrak t}}
\nc{\pw}{\widetilde{\mathfrak p}}
\nc{\qw}{\widetilde{\mathfrak q}}
\nc{\jw}{\widetilde j}

\nc{\grb}{\overline{\Gr}}
\nc{\I}{\mathcal I}

\nc{\lambdach}{{\check\lambda}}
\nc{\Lambdach}{{\check\Lambda}{}}
\nc{\much}{{\check\mu}}
\nc{\omegach}{{\check\omega}}
\nc{\nuch}{{\check\nu}}
\nc{\etach}{{\check\eta}}
\nc{\alphach}{{\check\alpha}}
\nc{\betach}{{\check\beta}}
\nc{\rhoch}{{\check\rho}}
\nc{\ch}{{\check h}}

\nc{\Hb}{\overline{\H}}

\nc{\BA}{{\mathbb{A}}}
\nc{\BC}{{\mathbb{C}}}
\nc{\BG}{{\mathbb{G}}}
\nc{\BM}{{\mathbb{M}}}
\nc{\BN}{{\mathbb{N}}}
\nc{\BP}{{\mathbb{P}}}
\nc{\BR}{{\mathbb{R}}}
\nc{\BZ}{{\mathbb{Z}}}
\nc{\BS}{{\mathbb{S}}}

\nc{\CA}{{\mathcal{A}}}
\nc{\CB}{{\mathcal{B}}}

\nc{\CK}{{\mathcal{K}}}

\nc{\CE}{{\mathcal{E}}}
\nc{\CF}{{\mathcal{F}}}
\nc{\CG}{{\mathcal{G}}}
\nc{\CI}{{\mathcal{I}}}
\nc{\CJ}{{\mathcal{J}}}
\nc{\CL}{{\mathcal{L}}}
\nc{\CM}{{\mathcal{M}}}
\nc{\CN}{{\mathcal{N}}}
\nc{\CO}{{\mathcal{O}}}
\nc{\CP}{{\mathcal{P}}}
\nc{\CQ}{{\mathcal{Q}}}
\nc{\CR}{{\mathcal{R}}}
\nc{\CS}{{\mathcal{S}}}
\nc{\CT}{{\mathcal{T}}}
\nc{\CU}{{\mathcal{U}}}
\nc{\CV}{{\mathcal{V}}}
\nc{\CW}{{\mathcal{W}}}
\nc{\CY}{{\mathcal Y}}
\nc{\CZ}{{\mathcal{Z}}}

\nc{\cM}{{\check{\mathcal M}}{}}
\nc{\csM}{{\check{\mathcal A}}{}}
\nc{\oM}{{\overset{\circ}{\mathcal M}}{}}
\nc{\obM}{{\overset{\circ}{\mathbf M}}{}}
\nc{\oCA}{{\overset{\circ}{\mathcal A}}{}}
\nc{\obA}{{\overset{\circ}{\mathbf A}}{}}
\nc{\ooM}{{\overset{\circ}{M}}{}}
\nc{\osM}{{\overset{\circ}{\mathsf M}}{}}
\nc{\vM}{{\overset{\bullet}{\mathcal M}}{}}
\nc{\nM}{{\underset{\bullet}{\mathcal M}}{}}
\nc{\oD}{{\overset{\circ}{\mathcal D}}{}}
\nc{\obD}{{\overset{\circ}{\mathbf D}}{}}
\nc{\oA}{{\overset{\circ}{\mathbb A}}{}}
\nc{\op}{{\overset{\bullet}{\mathbf p}}{}}
\nc{\cp}{{\overset{\circ}{\mathbf p}}{}}
\nc{\oU}{{\overset{\bullet}{\mathcal U}}{}}
\nc{\oZ}{{\overset{\circ}{\mathcal Z}}{}}
\nc{\ofZ}{{\overset{\circ}{\mathfrak Z}}{}}
\nc{\oF}{{\overset{\circ}{\fF}}}

\nc{\fa}{{\mathfrak{a}}}
\nc{\fb}{{\mathfrak{b}}}
\nc{\fg}{{\mathfrak{g}}}
\nc{\fgl}{{\mathfrak{gl}}}
\nc{\fh}{{\mathfrak{h}}}
\nc{\fj}{{\mathfrak{j}}}
\nc{\fk}{{\mathfrak{k}}}
\nc{\fm}{{\mathfrak{m}}}
\nc{\fn}{{\mathfrak{n}}}
\nc{\fu}{{\mathfrak{u}}}
\nc{\fp}{{\mathfrak{p}}}
\nc{\fr}{{\mathfrak{r}}}
\nc{\fs}{{\mathfrak{s}}}
\nc{\fsl}{{\mathfrak{sl}}}
\nc{\hsl}{{\widehat{\mathfrak{sl}}}}
\nc{\hgl}{{\widehat{\mathfrak{gl}}}}
\nc{\hg}{{\widehat{\mathfrak{g}}}}
\nc{\chg}{{\widehat{\mathfrak{g}}}{}^\vee}
\nc{\hn}{{\widehat{\mathfrak{n}}}}
\nc{\chn}{{\widehat{\mathfrak{n}}}{}^\vee}

\nc{\fA}{{\mathfrak{A}}}
\nc{\fB}{{\mathfrak{B}}}
\nc{\fD}{{\mathfrak{D}}}
\nc{\fE}{{\mathfrak{E}}}
\nc{\fF}{{\mathfrak{F}}}
\nc{\fG}{{\mathfrak{G}}}
\nc{\fK}{{\mathfrak{K}}}
\nc{\fL}{{\mathfrak{L}}}
\nc{\fM}{{\mathfrak{M}}}
\nc{\fN}{{\mathfrak{N}}}
\nc{\fP}{{\mathfrak{P}}}
\nc{\fU}{{\mathfrak{U}}}
\nc{\fV}{{\mathfrak{V}}}
\nc{\fZ}{{\mathfrak{Z}}}

\nc{\bb}{{\mathbf{b}}}
\nc{\bc}{{\mathbf{c}}}
\nc{\bd}{{\mathbf{d}}}
\nc{\be}{{\mathbf{e}}}
\nc{\bj}{{\mathbf{j}}}
\nc{\bn}{{\mathbf{n}}}
\nc{\bp}{{\mathbf{p}}}
\nc{\bq}{{\mathbf{q}}}
\nc{\bu}{{\mathbf{u}}}
\nc{\bv}{{\mathbf{v}}}
\nc{\bx}{{\mathbf{x}}}
\nc{\bs}{{\mathbf{s}}}
\nc{\by}{{\mathbf{y}}}
\nc{\bw}{{\mathbf{w}}}
\nc{\bA}{{\mathbf{A}}}
\nc{\bK}{{\mathbf{K}}}
\nc{\bB}{{\mathbf{B}}}
\nc{\bC}{{\mathbf{C}}}
\nc{\bD}{{\mathbf{D}}}
\nc{\bH}{{\mathbf{H}}}
\nc{\bM}{{\mathbf{M}}}
\nc{\bN}{{\mathbf{N}}}
\nc{\bV}{{\mathbf{V}}}
\nc{\bW}{{\mathbf{W}}}
\nc{\bX}{{\mathbf{X}}}
\nc{\bZ}{{\mathbf{Z}}}
\nc{\bS}{{\mathbf{S}}}

\nc{\sA}{{\mathsf{A}}}
\nc{\sB}{{\mathsf{B}}}
\nc{\sC}{{\mathsf{C}}}
\nc{\sD}{{\mathsf{D}}}
\nc{\sF}{{\mathsf{F}}}
\nc{\sK}{{\mathsf{K}}}
\nc{\sM}{{\mathsf{M}}}
\nc{\sO}{{\mathsf{O}}}
\nc{\sQ}{{\mathsf{Q}}}
\nc{\sP}{{\mathsf{P}}}
\nc{\sZ}{{\mathsf{Z}}}
\nc{\sfp}{{\mathsf{p}}}
\nc{\sr}{{\mathsf{r}}}
\nc{\sg}{{\mathsf{g}}}
\nc{\sff}{{\mathsf{f}}}
\nc{\sfb}{{\mathsf{b}}}
\nc{\sfc}{{\mathsf{c}}}
\nc{\sd}{{\mathsf{d}}}

\nc{\BK}{{\bar{K}}}

\nc{\tA}{{\widetilde{\mathbf{A}}}}
\nc{\tB}{{\widetilde{\mathcal{B}}}}
\nc{\tg}{{\widetilde{\mathfrak{g}}}}
\nc{\tG}{{\widetilde{G}}}
\nc{\TM}{{\widetilde{\mathbb{M}}}{}}
\nc{\tO}{{\widetilde{\mathsf{O}}}{}}
\nc{\tU}{{\widetilde{\mathfrak{U}}}{}}
\nc{\TZ}{{\tilde{Z}}}
\nc{\tx}{{\tilde{x}}}
\nc{\tbv}{{\tilde{\bv}}}
\nc{\tfP}{{\widetilde{\mathfrak{P}}}{}}
\nc{\tz}{{\tilde{\zeta}}}
\nc{\tmu}{{\tilde{\mu}}}

\nc{\urho}{\underline{\rho}}
\nc{\uB}{\underline{B}}
\nc{\uC}{{\underline{\mathbb{C}}}}
\nc{\ui}{\underline{i}}
\nc{\uj}{\underline{j}}
\nc{\ofP}{{\overline{\mathfrak{P}}}}
\nc{\oB}{{\overline{\mathcal{B}}}}
\nc{\og}{{\overline{\mathfrak{g}}}}
\nc{\oI}{{\overline{I}}}

\nc{\eps}{\varepsilon}
\nc{\hrho}{{\hat{\rho}}}

\nc{\one}{{\mathbf{1}}}
\nc{\two}{{\mathbf{t}}}

\nc{\Rep}{{\mathop{\operatorname{\rm Rep}}}}
\nc{\Tot}{{\mathop{\operatorname{\rm Tot}}}}
\nc{\Ker}{{\mathop{\operatorname{\rm Ker}}}}
\nc{\Hilb}{{\mathop{\operatorname{\rm Hilb}}}}
\nc{\End}{{\mathop{\operatorname{\rm End}}}}
\nc{\Ext}{{\mathop{\operatorname{\rm Ext}}}}
\nc{\CHom}{{\mathop{\operatorname{{\mathcal{H}}\it om}}}}
\nc{\GL}{{\mathop{\operatorname{\rm GL}}}}
\nc{\gr}{{\mathop{\operatorname{\rm gr}}}}
\nc{\Id}{{\mathop{\operatorname{\rm Id}}}}
\nc{\de}{{\mathop{\operatorname{\rm def}}}}
\nc{\length}{{\mathop{\operatorname{\rm length}}}}
\nc{\supp}{{\mathop{\operatorname{\rm supp}}}}

\nc{\Cliff}{{\mathsf{Cliff}}}
\nc{\Fl}{\on{Fl}}
\nc{\Fib}{{\mathsf{Fib}}}
\nc{\Coh}{{\mathsf{Coh}}}
\nc{\FCoh}{{\mathsf{FCoh}}}

\nc{\reg}{{\text{\rm reg}}}

\nc{\cplus}{{\mathbf{C}_+}}
\nc{\cminus}{{\mathbf{C}_-}}
\nc{\cthree}{{\mathbf{C}_*}}
\nc{\Qbar}{{\bar{Q}}}

\nc{\bh}{{\bar{h}}}
\nc{\bOmega}{{\overline{\Omega}}}

\nc{\seq}[1]{\stackrel{#1}{\sim}}

\nc{\chA}{\check A}
\nc{\chrho}{\check \rho}
\nc{\chT}{{\check T}}
\nc{\chM}{{\check M}}
\nc{\chG}{{\check G}}
\nc{\chH}{{\check H}}

\nc{\Stab}{{\mathop{\operatorname{\rm Stab }}}}

\nc{\domwts}{{\check\Lambda^+}}
\nc{\tdomwts}{\tilde{\check{\Lambda}}^+}
\nc{\wts}{{\check\Lambda}}
\nc{\twts}{{\tilde{\check{\Lambda}}}}

\nc{\domcowts}{\Lambda^+}
\nc{\cowts}{\Lambda}
\nc{\poscowts}{\Lambda^{\on{pos}}}
\nc{\poscoroots}{R^{\on{pos}}}

\nc{\lat}{\Lambda}
\nc{\poslat}{\Lambda^{\on{pos}}}

\nc{\blambda}{{\bar\lambda}}
\nc{\risom}{\stackrel{\sim}{\rightarrow}}
\nc{\lisom}{\stackrel{\sim}{\leftarrow}}
\nc{\hAX}{{{\hat A}_X}}
\nc{\Zhecke}{\hat{Z}}
\nc{\Vect}{\on{Vect}}

\nc{\can}{{\on{can}}}

\nc{\Z}{{Z}}
\nc{\stZ}{{^\star Z}}
\nc{\nuZ}{{{}_{c,\geq \nu}Z}}
\nc{\zZ}{{{}_{c,\geq 0}Z}}

\nc{\Zcan}{{Z_{\can}}}
\nc{\bZcan}{{\overline{Z}}_{\can}}
\nc{\tZcan}{{\widetilde{Z}}_{\can}}
\nc{\stZcan}{{{}^\star Z_{\can}}}

\nc{\bBunP}{{{}_{\infty}\overline{\Bun}}_P}
\nc{\bzBunP}{{{}_{\infty,\geq 0}\overline{\Bun}}_P}
\nc{\bnuBunP}{{{}_{\infty,\geq \nu}\overline{\Bun}}_P}
\nc{\tBunP}{{{}_{\infty}\widetilde{\Bun}}_P}
\nc{\tnuBunP}{{{}_{\infty,\geq \nu}\widetilde{\Bun}}_P}
\nc{\tzBunP}{{{}_{\infty,\geq 0}\widetilde{\Bun}}_P}

\nc{\bBunB}{{{}_{\infty}\overline{\Bun}}_B}

\nc{\hrG}{h^{\rightarrow}_G}
\nc{\hlG}{h^{\leftarrow}_G}
\nc{\cHG}{{{}_c\H_G}}
\nc{\cHH}{{{}_c\H_H}}

\nc{\gtimes}{\stackrel{G}{\times}}
\nc{\btimes}{\stackrel{B}{\times}}
\nc{\timesBunG}{\underset{\Bun_G}{\times}}
\nc{\timesBunM}{\underset{\Bun_M}{\times}}
\nc{\timesBunP}{\underset{\Bun_P}{\times}}

\nc{\ctheta}{{{}_c\theta}}
\nc{\ceta}{{{}_c\eta}}
\nc{\ckappa}{{{}_c\kappa}}

\nc{\sJ}{{\mathsf{J}}}
\nc{\Junk}{{\mathsf{Junk}}}

\nc{\catp}{{\on{P}}}
\nc{\catq}{{\on{Q}}}
\nc{\openX}{{\mathring X}}

\nc{\disj}{{\on{disj}}}

\nc{\hl}{h^{\leftarrow}}
\nc{\hr}{h^{\rightarrow}}

\nc{\affgr}{\on{Gr}}

\begin{document}

\title[Horospherical varieties]
{Hecke operators on quasimaps\\ into horospherical varieties}

\author[Gaitsgory and Nadler]{Dennis Gaitsgory and David Nadler}
\address{Harvard University, Northwestern University}
\email{gaitsgde@math.harvard.edu, nadler@math.northwestern.edu}

\begin{abstract}
Let $G$ be a connected reductive complex algebraic group.
This paper and its companion~\cite{GNcombo06} are
devoted to the space $Z$ of meromorphic quasimaps
from a curve into an affine spherical $G$-variety $X$. The space $Z$ may be thought
of as an algebraic
model for the loop space of $X$.
The theory we develop
associates to $X$ a connected reductive complex algebraic subgroup $\check H$
of the dual group $\check G$. The construction of $\check H$ is via
Tannakian formalism: we identify
a certain tensor category $\catq(Z)$ of perverse sheaves on $Z$ with the category of finite-dimensional representations of $\check H$.

In this paper, we focus on horospherical varieties, a class of varieties closely related to flag varieties. 
For an affine horospherical $G$-variety $X_{\on{horo}}$,
the category $\catq(Z_{\on{horo}})$ is equivalent to
a category of vector spaces
graded by a lattice. Thus the associated subgroup $\check H_{\on{horo}}$ is a torus.
%
The case of horospherical varieties may be thought of as a simple example, but it also
plays a central role in the general theory. To an arbitrary affine spherical $G$-variety $X$,
one may associate a horospherical variety $X_{\on{horo}}$.
Its associated subgroup $\check H_{\on{horo}}$
turns out to be a maximal torus in the subgroup $\check H$ associated to $X$.

\end{abstract}

\maketitle



\begin{section}{Introduction}

Let $G$ be a connected reductive complex algebraic group.
In this paper and its companion~\cite{GNcombo06},
we study the space $Z$ of meromorphic quasimaps
from a curve into an affine spherical $G$-variety $X$. 
A $G$-variety $X$ is said to be spherical  if a Borel subgroup of $G$
acts on $X$ with a dense orbit. Examples include flag varieties,
symmetric spaces, and toric varieties.
A meromorphic quasimap consists of a point of the curve, a $G$-bundle on the curve,
and a meromorphic section of 
the associated $X$-bundle with a pole only at the distinguished point.
The space $Z$ may be thought
of as an algebraic
model for the loop space of $X$.

The theory we develop identifies a certain tensor category $\catq(Z)$ 
of perverse sheaves on $Z$
with the category of finite-dimensional representations
of a connected reductive complex algebraic
subgroup $\check H$ of the dual group $\check G$.
Our method is to use Tannakian formalism:
we endow $\catq(Z)$ with a tensor product, a fiber functor to vector spaces,
and the necessary compatibility constraints so that it must be equivalent to
the category of representations of such a group.
Under this equivalence, the fiber functor corresponds
to the forgetful functor 
which assigns to a representation of $\check H$ its underlying vector space.
In the paper~\cite{GNcombo06}, we define the category $\catq(Z)$,
and endow it with a tensor product and fiber functor.
This paper provides a key technical result needed for the construction
of the fiber functor.

Horospherical $G$-varieties form a special class of $G$-varieties closely related to flag varieties.
A subgroup $S\subset G$
is said to be horospherical if it contains the unipotent radical of a Borel subgroup of $G$.
A $G$-variety $X$ is said to be horospherical if for
each point $x\in X$, its stabilizer $S_x\subset G$ is horospherical.
When $X$ is an affine horospherical $G$-variety,
the subgroup $\check H$ we associate to it turns out to be a torus.
To see this, we explicitly calculate the functor which corresponds to
the restriction of representations from $\check G$.
Representations of $\check G$ naturally act on the category $\catq(Z)$
via the geometric Satake correspondence.
The restriction of representations is given by applying this action to the object
of $\catq(Z)$ corresponding to the trivial representation of $\check H$.
The main result of this paper
describes this action in the horospherical case.
The statement does not mention $\catq(Z)$, 
but rather what is needed in~\cite{GNcombo06}
where we define and study $\catq(Z)$.

In the remainder of the introduction, we first describe a piece of the theory
of geometric Eisenstein series which the main result of this paper generalizes.
This may give the reader some context from which to approach 
the space $Z$ and our main result.
We then define $Z$ and state our main result.
Finally, we collect notation and preliminary results needed in what follows.
Throughout the introduction, we use the term space 
for objects which are strictly speaking stacks and ind-stacks.


\begin{subsection}{Background}
One way to approach the results of this paper
is to interpret them as a
generalization of a theorem of Braverman-Gaitsgory~\cite[Theorem 3.1.4]{BG02}
from the theory of geometric Eisenstein series.
Let $C$ be a smooth complete complex algebraic curve.
The primary aim of the geometric Langlands program is to
construct sheaves on the moduli space $\Bun_G$ of $G$-bundles on $C$
which are eigensheaves
for Hecke operators. 
These 
are the operators which result from modifying $G$-bundles
at prescribed points of the curve $C$.
Roughly speaking, the theory of geometric Eisenstein series constructs sheaves on $\Bun_G$
starting with local systems
on the moduli space $\Bun_T$,
where $T$ is the universal Cartan of $G$.
When the original local system is sufficiently generic, 
the resulting sheaf is an eigensheaf for the Hecke operators.

At first glance,
the link between $\Bun_T$ and $\Bun_G$ should be the moduli stack $\Bun_B$ of $B$-bundles
on $C$, where $B\subset G$ is a Borel subgroup with unipotent radical $U\subset B$
and reductive quotient $T=B/U$. Unfortunately,
naively working with the natural diagram
$$
\begin{array}{ccc}
\Bun_B & \to & \Bun_G \\
\downarrow & & \\
\Bun_T & & \\
\end{array}
$$
leads to difficulties: the fibers of the horizontal map are not compact.
The eventual successful construction depends on V.~Drinfeld's relative
compactification of $\Bun_B$ along the fibers of the map to $\Bun_G$.
The starting point for the compactification is the
observation that $\Bun_B$ also classifies data 
$$
(\CP_G\in\Bun_G,\CP_T\in\Bun_T,\sigma:
\CP_{T}\to \CP_G{\gtimes}{G/U})
$$
where $\sigma$ is a $T$-equivariant bundle map
to the $\CP_G$-twist of $G/U$.
From this perspective, it is natural to be less restrictive and allow maps into the $\CP_G$-twist of
the fundamental affine space 
$$
\ol{G/U}=\Spec(\BC[G]^U).
$$
Here $\BC[G]$ denotes the ring of regular functions on $G$,
and $\BC[G]^U\subset\BC[G]$ the (right) $U$-invariants.
Following V. Drinfeld, we define the compactification $\ol{\Bun}_B$ to be that classifying
quasimaps
$$
(\CP_G\in \Bun_G,\CP_T\in\Bun_T,\sigma:\CP_{T}\to 
\CP_G{\gtimes}\overline{G/U})
$$
where 
$\sigma$ is a $T$-equivariant bundle map which factors 
$$
\sigma|_{C'}:\CP_{T}|_{C'}\to\CP_G{\gtimes}{G/U}|_{C'}\to 
\CP_G{\gtimes}\overline{G/U}|_{C'},
$$
for some open curve $C'\subset C$.
Of course, 
the quasimaps that satisfy
$$
\sigma:\CP_{T}\to 
\CP_G{\gtimes}{G/U}
$$
form a subspace canonically isomorphic to $\Bun_B$.

Since the Hecke operators on $\Bun_G$ do not lift to $\ol{\Bun}_B$,
it is useful to introduce a version of $\ol{\Bun}_B$ on which they do.
Following \cite[Section 4]{BG02}, 
we define the space $\bBunB$ to be
that classifying meromorphic quasimaps
$$
(c\in C,\CP_G\in \Bun_G,\CP_T\in\Bun_T,\sigma:\CP_{T}|_{C\setminus c}\to 
\CP_G{\gtimes}\overline{G/U}|_{C\setminus c})
$$
where 
$\sigma$ is a $T$-equivariant bundle map which factors 
$$
\sigma|_{C'}:\CP_{T}|_{C'}\to\CP_G{\gtimes}{G/U}|_{C'}\to 
\CP_G{\gtimes}\overline{G/U}|_{C'},
$$
for some open curve $C'\subset C\setminus c$.
We call $c\in C$ the pole point of the quasimap.
Given a meromorphic quasimap with $G$-bundle $\CP_G$ and pole point $c\in C$, 
we may modify $\CP_G$ at $c$ and obtain a new meromorphic quasimap.
In this way, the Hecke operators on $\Bun_G$ lift to $\bBunB$.

Now the result we seek to generalize \cite[Theorem 3.1.4]{BG02} 
describes how
the Hecke operators act on a distinguished object of the 
category $\catp(\bBunB)$ of perverse sheaves with $\BC$-coefficients
on $\bBunB$.
Let $\cowts_G=\Hom(\BC^\times,T)$ be the coweight lattice,
and let $\domcowts_G\subset\cowts$ be the semigroup of dominant coweights of $G$.
For $\lambda\in\domcowts_G$, we have the Hecke operator
$$
H^\lambda_G:\catp(\bBunB)
\to\catp(\bBunB)
$$ 
given by convolving with the simple spherical modification of coweight $\lambda$. 
(See \cite[Section 4]{BG02} or
Section~\ref{sconv} below for more details.)
For $\mu\in\cowts_G$,
we have the locally closed subspace
$\bBunB^{\mu}\subset \bBunB$ that classifies data
for which the map 
$$
\CP_{T}(\mu\cdot c)|_{C\setminus c}
\stackrel{\sigma}{\to} 
\CP_G{\gtimes}\overline{G/U}|_{C\setminus c}
$$
extends to a holomorphic map 
$$
\CP_{T}(\mu\cdot c)
\stackrel{\sigma}{\to} 
\CP_G{\gtimes}\overline{G/U} 
$$
which factors
$$
\CP_{T}(\mu\cdot c)
\stackrel{\sigma}{\to} 
\CP_G{\gtimes}{G/U}\to
\CP_G{\gtimes}\overline{G/U}.
$$
We write 
$\bBunB^{\leq\mu}\subset\bBunB$ for the closure of 
$\bBunB^{\mu}\subset\bBunB$,
and 
$$
\IC^{\leq\mu}_{\bBunB}\in\catp(\bBunB)
$$ 
for the intersection cohomology sheaf of $\bBunB^{\leq\mu}\subset\bBunB$.

\begin{thm}\cite[Theorem 3.1.4]{BG02}
For $\lambda\in\domcowts_G$, there is a canonical isomorphism
$$
H^\lambda_G (\IC^{\leq 0}_{\bBunB})\simeq
\sum_{\mu\in\cowts_T}
\IC^{\leq\mu}_{\bBunB} \otimes \Hom_{\chT} (V_\chT^\mu,V_\chG^\lambda)
$$

\end{thm}

Here we write $V^\lambda_\chG$ for the irreducible
representation of the dual group $\chG$ with highest weight $\lambda\in\domcowts_G$,
and $V^\mu_\chT$ for the irreducible
representation of the dual torus $\chT$ of weight $\mu\in\cowts_G$.

In the same paper of Braverman-Gaitsgory~\cite[Section 4]{BG02}, 
there is a generalization \cite[Theorem 4.1.5]{BG02}
of this theorem from the Borel
subgroup $B\subset G$ to
other parabolic subgroups $P\subset G$.
We recall and use this generalization in Section~\ref{sconv} below.
It is the starting point for the results of this paper.

\end{subsection}


\begin{subsection}{Main result}
The main result of this paper is a version of \cite[Theorem 3.1.4]{BG02} for $X$ an
arbitrary affine horospherical
$G$-variety with a dense $G$-orbit $\openX\subset X$.
For any point in the dense $G$-orbit $\openX\subset X$, 
we refer to its stabilizer $S\subset G$ as the generic stabilizer of $X$. 
All such subgroups are conjugate to each other.
By choosing such a point, we obtain an identification
$
\openX\simeq G/S.
$

To state our main theorem, we first introduce some more notation.
Satz 2.1 of \cite{Kinv90}
states that the normalizer of a horospherical subgroup $S\subset G$
is a parabolic subgroup $P\subset G$
with the same derived group $[P,P]=[S,S]$.
We write $A$ for the quotient torus $P/S$,
and $\cowts_A=\Hom(\BC^\times, A)$ for its coweight lattice. 
Similarly, for the identity component $S^0\subset S$,
we write $A_0$ for the quotient torus $P/S_0$,
and $\cowts_{A_0}=\Hom(\BC^\times, A_0)$ for its coweight lattice. 
The natural maps $T\to A_0\to A$ induce maps of coweight lattices
$$
\cowts_T\stackrel{q}{\to}\cowts_{A_0}\stackrel{i}{\to}\cowts_A,
$$
where $q$ is a surjection, and $i$ is an injection.
For a conjugate of $S$, the associated tori are canonically isomorphic
to those associated to $S$.
Thus when $S$ is the generic stabilizer of a horospherical $G$-variety $X$,
the above tori, lattices and maps are canonically associated to $X$.

For an affine horospherical $G$-variety $X$ with dense $G$-orbit $\openX\subset X$, 
we define the space $Z$ to be that classifying mermorphic quasimaps into $X$.
Such a quasimap consists of data
$$
(c\in C,\CP_G\in\Bun_G,\sigma:{C\setminus c}\to 
\CP_G{\gtimes}X|_{C\setminus c})
$$
where 
$\sigma$ is a section which
factors
$$
\sigma|_{C'}:C'\to 
\CP_G{\gtimes}\openX|_{C'}\to
\CP_G{\gtimes}{X}|_{C'},
$$
for some open curve $C'\subset C\setminus c$.

Given a meromorphic quasimap into $X$ with $G$-bundle $\CP_G$ and pole point $c\in C$, 
we may modify $\CP_G$ at $c$ and obtain a new meromorphic quasimap.
But in this context the resulting Hecke operators on $Z$ do not in general preserve the category 
of perverse sheaves. Instead, we must consider the bounded derived category $\Sh(Z)$ of sheaves of $\BC$-modules
on $Z$.
For $\lambda\in\domcowts_G$, we have the Hecke operator
$$
H^\lambda_G:\Sh(Z)
\to\Sh(Z)
$$ 
given by convolving with the simple spherical modification of coweight $\lambda$.
(See Section~\ref{sconv} below for more details.)
For $\kappa\in\cowts_{A_0}$, we have a 
locally closed subspace $Z^{\kappa}\subset Z$ consisting of meromorphic quasimaps
that factor
$$
\sigma:C\setminus c\to 
\CP_G{\gtimes}\openX|_{C\setminus c}\to
\CP_G{\gtimes}{X}|_{C\setminus c}
$$
and have a singularity of type $\kappa$ at $c\in C$.
(See Section~\ref{sectnaive} below for more details.)
We write $Z^{\leq\kappa}\subset Z$
for the closure of $Z^{\kappa}\subset Z$,
and 
$$
\IC_Z^{\leq\kappa}\in\Sh(Z)
$$ 
for its intersection cohomology sheaf. 

Our main result is the following.

\begin{thm}
For $\lambda\in\domcowts_G$, there is an isomorphism
$$
H^\lambda_G (\IC_{Z}^{\leq 0})\simeq
\sum_{\kappa\in\cowts_{A_0}}
\sum_{
{\mu\in\cowts_{T}},
{q(\mu)=\kappa}
}
\IC_{Z}^{\leq\kappa} \otimes \Hom_{\chT} (V_\chT^\mu,V_\chG^\lambda)
[\langle2\chrho_M,\mu\rangle].
$$
\end{thm}

Here the torus $A_0$ and its coweight lattice $\cowts_{A_0}$ are those 
associated to the generic stabilizer $S\subset G$.
We write $M$ for the Levi quotient of the normalizer $P\subset G$
of the generic stabilizer $S\subset G$, and
$2\check\rho_M$ for the sum of the positive roots of $M$.

\end{subsection}


\begin{subsection}{Notation}

Throughout this paper, let $G$ be a connected reductive complex algebraic group,
let $B\subset G$ be a  Borel subgroup with unipotent radical $U(B)$, and 
let $T=B/U(B)$ be the abstract Cartan. 

Let $\wts_G$ denote the weight lattice $\Hom(T,\BC^\times)$, and
$\domwts_G\subset\wts_G$ the semigroup of dominant weights.
For $\lambda\in\domwts_G$,
we write $V_{G}^\lambda$ for the irreducible
representation of $G$ of highest weight $\lambda$.

Let $\cowts_G$ denote the coweight lattice $\Hom(\BC^\times, T)$, and
$\domcowts_G\subset \cowts_G$ the semigroup of dominant coweights.
For $\lambda\in\domcowts_G$,
let $V_{\chG}^\lambda$ denote the irreducible
representation of the dual group $\chG$ of highest weight $\lambda$.

Let $\poscowts_G\subset\cowts_G$ denote the semigroup of coweights in $\cowts_G$
which are non-negative on $\domwts_G$, and let $\poscoroots_G\subset\poscowts_G$ 
denote the semigroup of positive coroots.

Let $P\subset G$ be a parabolic subgroup
with unipotent radical $U(P)$,
and let $M$ be the Levi factor $P/U(P)$. 

We have the natural map 
$$
\check r:\wts_{M/[M,M]}\to\wts_G
$$ of weights,
and the dual map 
$$
r:\cowts_G\to\cowts_{M/[M,M]}
$$ of coweights.

Let $\domwts_{G,P}\subset\wts_{M/[M,M]}$ 
denote the inverse image $\check r^{-1}(\domwts_G)$.
Let $\poscowts_{G,P}\subset\cowts_{M/[M,M]}$ 
denote the semigroup of coweights in $\cowts_{M/[M,M]}$
which are non-negative on $\domwts_{G,P}$.
Let $\poscoroots_{G,P}\subset \poscowts_{G,P}$ denote the image $r(\poscoroots_G)$.

Let $\CW_M$ denote the Weyl group of $M$, and let 
$\CW_M\domwts_G\subset\wts_G$ denote the union of the translates
of $\domwts_G$ by $\CW_M$.
Let
$\tilde{\Lambda}^\pos_{G,P}\subset\domcowts_M$ denote
the semigroup of dominant coweights of $M$ which are nonnegative
on $\CW_M\domwts_G$.

Finally, let $\langle\cdot,\cdot\rangle:\wts_G\times\cowts_G\to\BZ$ denote the natural pairing,
and let $\chrho_M\in\wts_G$ denote half the sum of the positive roots of $M$.

\end{subsection}


\begin{subsection}{Bundles and Hecke correspondences}\label{sbandh}

Let $C$ be a smooth complete complex algebraic curve.

For a connected complex algebraic group $H$,
let $\Bun_H$ be the moduli stack
of $H$-bundles on $C$. 
Objects of $\Bun_H$ will be denoted by $\CP_H$.

Let $\H_H$ be the Hecke ind-stack that classifies data
$$
(c\in C,\CP_H^1,\CP^2_H\in\Bun_H,
\alpha:\CP^1_H|_{C\setminus c}\stackrel{\sim}{\to} 
\CP^2_H|_{C\setminus c})
$$
where $\alpha$ is an isomorphism of $H$-bundles.
We have the maps 
$$
\Bun_H\stackrel{\hl_H}{\leftarrow}\H_H\stackrel{\hr_H}{\to}\Bun_H
$$
defined by 
$$
\hl_H(c,\CP_H^1,\CP^2_H,\alpha)=\CP^1_H\qquad
\hr_H(c,\CP_H^1,\CP^2_H,\alpha)=\CP^2_H,
$$
and the map
$$
\pi:\H_H\to C
$$
defined by
$$
\pi(c,\CP_H^1,\CP^2_H,\alpha)=c.
$$

It is useful to have another description of the Hecke ind-stack $\H_H$
for which we introduce some more notation.
Let $\CO$ be the ring of formal power series $\BC[[t]]$, let $\CK$
be the field of formal Laurent series $\BC((t))$, and let $D$ be the formal
disk $\Spec(\CO)$. 
For a point $c\in C$, let $\CO_c$ be the completed local ring
of $C$ at $c$, 
and let $D_c$ be the formal
disk $\Spec(\CO_c)$. 
Let $\Aut(\CO)$ be the group-scheme of automorphisms of the ring $\CO$.
Let $H(\CO)$ be the group of $\CO$-valued points of $H$,
and let $H(\CK)$ be the group of $\CK$-valued points of $H$.
Let $\affgr_{H}$ be the affine Grassmannian of $H$. It is an ind-scheme
whose set of $\BC$-points is the quotient $H(\CK)/H(\CO)$.

Now consider the $(H(\CO)\rtimes\Aut(\CO))$-torsor
$$
\widehat{\Bun_H\times C}\to\Bun_H\times C
$$
that classifies data
$$
(c\in C,\CP_H\in\Bun_H,
\beta:D\times H\risom\CP_H|_{D_c},\gamma:D\risom D_c)
$$
where $\beta$ is an isomorphism of $H$-bundles,
and $\gamma$ is an identification of formal disks.
We have an identification
$$
\H_H \simeq \widehat{\Bun_H\times C}
\overset{(H(\CO)\rtimes\Aut(\CO))}{\times}
\affgr_{H}
$$
such that the projection $\hr_H$ corresponds to the obvious
projection from the twisted product to $\Bun_H$.

For $H$ reductive, the $(H(\CO)\rtimes\Aut(\CO))$-orbits $\affgr_H^\lambda\subset\affgr_H$
are indexed by $\lambda\in\domcowts_H$.
For $\lambda\in\domcowts_H$,
we write $\H_H^\lambda\subset \H_H$ for the substack
$$
\H_H^\lambda \simeq \widehat{\Bun_H\times C}
\overset{(H(\CO)\rtimes\Aut(\CO))}{\times}
\affgr^\lambda_{H}.
$$

For a parabolic subgroup $P\subset H$, 
the connected components $S_{P,\theta}\subset\affgr_P$ are indexed by
$\theta\in \cowts_{P}/\cowts_{[P,P]^{sc}}$, where $[P,P]^{sc}$ denotes the simply connected cover
of $[P,P]$.
For $\theta\in\cowts_{P}/\cowts_{[P,P]^{sc}}$,
we write $\CS_{P,\theta}\subset \H_P$ for the ind-substack
$$
\CS_{P,\theta} \simeq \widehat{\Bun_P\times C}
\overset{(P(\CO)\rtimes\Aut(\CO))}{\times}
S_{P,\theta}.
$$
For $\theta\in\cowts_{P}/\cowts_{[P,P]^{sc}}$, and $\lambda\in\domcowts_H$,
we write $\CS^\lambda_{P,\theta}\subset \H_P$ for the ind-substack
$$
\CS^\lambda_{P,\theta} \simeq \widehat{\Bun_P\times C}
\overset{(P(\CO)\rtimes\Aut(\CO))}{\times}
S^\lambda_{P,\theta}
$$
where $S_{P,\theta}^\lambda$ denotes the intersection $S_{P,\theta}\cap\affgr_H^\lambda$.

For any ind-stack $\CZ$ over $\Bun_H\times C$,
we have the $(H(\CO)\rtimes\Aut(\CO))$-torsor
$$
\widehat{\CZ}\to\CZ
$$
obtained by pulling back 
the $(H(\CO)\rtimes\Aut(\CO))$-torsor
$$
\widehat{\Bun_H\times C}\to\Bun_H\times C.
$$
We also have the Cartesian diagram
$$
\begin{array}{ccc}
\H_H\underset{\Bun_H\times C}{\times}\CZ
& \stackrel{\hr_H}{\to} & \CZ \\
\downarrow & & \downarrow \\
\H_H & \stackrel{\hr_H}{\to} & \Bun_H
\end{array}
$$
and
an identification
$$
\H_H\underset{\Bun_H\times C}{\times}\CZ\simeq
\widehat{\CZ}
\overset{(H(\CO)\rtimes\Aut(\CO))}{\times}
\affgr_{H}
$$
such that the projection $\hr_H$ corresponds to the obvious
projection from the twisted product to $\CZ$.
For 
$\CF\in\Sh(\CZ)$, and 
$\CP\in\catp_{(H(\CO)\rtimes\Aut(\CO))}(\affgr_H)$,
we may form the twisted product
$$
(\CF\tboxtimes\CP)^r\in\Sh(\H_H\underset{\Bun_H\times C}{\times}\CZ).
$$
with respect to the map $\hr_H$.
In particular, for $\lambda\in\domcowts_H$, we may take $\CP$ to be 
the intersection cohomology sheaf $\CA_G^\lambda$ of the closure 
$\ol\affgr_H^\lambda\subset\affgr_H$ of the $(H(\CO)\rtimes\Aut(\CO))$-orbit
$\affgr_H^\lambda\subset\affgr_H$.

\end{subsection}

\end{section}


\begin{section}{Affine horospherical $G$-varieties}\label{sectaffvar}

A subgroup $S\subset G$ is said to be horospherical
if it contains the unipotent radical of a Borel subgroup of $G$.
A $G$-variety $X$ is said to be {horospherical} if for each point $x\in X$,
its stabilizer $S_x\subset G$ is horospherical.
A $G$-variety $X$ is said to be {spherical} if a Borel subgroup of $G$
acts on $X$ with a dense orbit.
Note that a horospherical $G$-variety contains a dense $G$-orbit
if and only if it is spherical.

Let $X$ be an affine $G$-variety.
As a representation of $G$, the ring of regular functions $\BC[X]$
decomposes into isotypic components
$$
 \BC[X]\simeq \sum_{\lambda\in \domwts_G} \BC[X]_\lambda.
$$
We say that $\BC[X]$ is {graded} if
$$
\BC[X]_\lambda\BC[X]_\mu\subset \BC[X]_{\lambda+\mu},
$$
for all $\lambda,\mu\in\domwts_G$.
We say that $\BC[X]$ is {simple} if the irreducible representation
$V^\lambda$ of highest weight $\lambda$
occurs in $\BC[X]_\lambda$ with multiplicity $0$ or $1$, for
all $\lambda\in\domwts_G$.

\begin{prop} Let $X$ be an affine $G$-variety.

(1) \cite[Proposition 8, (3)]{Pmsb86}
 $X$ is horospherical if and only if $\BC[X]$ is graded.

(2) \cite[Theorem 1]{Pmsb86}  
$X$ is spherical if and only if $\BC[X]$ is simple.

\end{prop}

We see by the proposition that affine horospherical $G$-varieties containing a dense $G$-orbit are classified by
finitely-generated subsemigroups of $\domwts_G$.
To such a variety $X$, one associates the subsemigroup 
$$
\domwts_X\subset\domwts_G
$$ 
of dominant weights $\lambda$ with $\dim \BC[X]_\lambda >0$.


\begin{subsection}{Structure of generic stabilizer}

\begin{thm}\cite[Satz 2.2]{Kinv90}\label{thmopenset}
If $X$ is an irreducible horospherical $G$-variety,
then there is an open $G$-invariant subset $W\subset X$, 
and a $G$-equivariant isomorphism
$W\simeq G/S\times Y$, where $S\subset G$ is a horospherical subgroup,
and $Y$ is a variety on which $G$ acts trivially.
\end{thm}

Note that for any two such open subsets $W\subset X$ and isomorphisms
$W\simeq G/S\times Y$, the subgroups $S\subset G$ are conjugate.
We refer to such a subgroup $S\subset G$ as the generic stabilizer of $X$.

\begin{lem}\label{lemstab}\cite[Satz 2.1]{Kinv90}
If $S\subset G$ is a horospherical subgroup,
then its normalizer is a parabolic subgroup $P\subset G$ with the same
derived group $[P,P]=[S,S]$ and unipotent radical $U(P)=U(S)$.
\end{lem}

Note that the identity component $S^0\subset S$ is also horospherical
with the same derived group
$[S^0,S^0]=[S,S]$ and unipotent radical $U(S^0)=U(S)$.

Let $S\subset G$ be a horospherical subgroup with identity component $S^0\subset S$,
and normalizer $P\subset G$.
We write $A$ for the quotient torus $P/S$,
and $\cowts_A$ for its coweight lattice $\Hom(\BC^\times, A)$. 
Similarly, 
we write $A_0$ for the quotient torus $P/S^0$,
and $\cowts_{A_0}$ for its coweight lattice $\Hom(\BC^\times, A_0)$. 
The natural maps 
$$T\to A_0\to A
$$ 
induce maps of coweight lattices
$$
\cowts_T\stackrel{q}{\to}\cowts_{A_0}\stackrel{i}{\to}\cowts_A,
$$
where $q$ is a surjection, and $i$ is an injection.
For a conjugate of $S$, the associated tori, lattices, and maps are canonically isomorphic 
to those associated to $S$.
Thus when $S$ is the generic stabilizer of a horospherical $G$-variety $X$,
the tori, lattices and maps are canonically associated to $X$.

We shall need the following finer description of which subgroups $S\subset G$
may appear as the generic stabilizer of an affine horospherical $G$-variety.
To state it, we introduce some more notation used throughout the paper.
For a horospherical subgroup $S\subset G$ with identity component $S^0\subset S$,
and normalizer $P\subset G$,
let $M$ be the Levi quotient $P/U(P)$, 
let $M_S$ be the Levi quotient $S/U(S)$, 
and let $M_S^0$ be the identity component
of $M_S$. 
The natural maps
$$
S^0\to S\to P
$$
induce isomorphisms of
derived groups 
$$
[M^0_S,M_S^0]\risom [M_S,M_S]\risom [M,M].
$$ 
We write $\cowts_{M/[M,M]}$ for the coweight lattice of the torus $M/[M,M]$,
and $\cowts_{M^0_S/[M_S,M_S]}$ for the coweight lattice of the torus
$M^0_S/[M_S,M_S]$.
The natural maps 
$$
M^0_S/[M_S,M_S]\to M/[M,M]\to A_0
$$
induce a short exact sequence of coweight lattices
$$
0\to\cowts_{M^0_S/[M_S,M_S]}
\to
\cowts_{M/[M,M]}
\to
\cowts_{A_0}
\to 0.
$$

%

\begin{prop}\label{propstab}
Let $S\subset G$ be a horospherical subgroup.
Then $S$ is the generic stabilizer of an affine horospherical $G$-variety
containing a dense $G$-orbit
if and only if 
$$
\cowts_{M^0_S/[M_S,M_S]}\cap\poscowts_{G,P}=\langle0\rangle.
$$
\end{prop}

\begin{proof}
The proof of the proposition relies on the following lemma.
Let $\check V$ be a finite-dimensional real
vector space, and let $\check V^+$ be an open 
set in $\check V$
which is preserved by the action of $\BR^{>0}$.
Let $V$ be the dual of $\check V$, and let $V^\pos$ be the closed cone of covectors
in $V$ that are nonnegative on all vectors in $\check V^+$.
For a linear subspace $\check W\subset \check V$, we write $\check W^\perp\subset V$
for its orthogonal.

\begin{lem}\label{lemortho}
The map $\check W\mapsto \check W^\perp$ provides a bijection
from the set of all linear subspaces $\check W\subset \check V$ 
such that $\check W\cap \check V^+\neq\emptyset$
to the set of all linear subspaces
$W\subset V$ such that $W\cap V^\pos=\langle0\rangle$.
\end{lem}

\begin{proof}
If $\check W\cap\check V^+\neq\emptyset$, then clearly
$\check W^\perp\cap V^\pos=\langle0\rangle$.
Conversely, if $W\cap V^\pos=\langle0\rangle$, then since $\check V^+$ is open,
there is a hyperplane $H\subset V$
such that $W\subset H$, and $H\cap V^\pos=\langle0\rangle$. 
Thus $H^\perp\subset W^\perp$, and $H^\perp\cap \check V^+\neq\emptyset$,
and so $W^\perp \cap \check V^+\neq\emptyset$.
\end{proof}

Now suppose $X$ is an affine horospherical $G$-variety with an open $G$-orbit
and generic stabilizer 
$S\subset G$ with normalizer $P\subset G$.
Then we have
$\domwts_X\subset \domwts_{G,P}$, since otherwise $[S,S]$ would be smaller.
We also have that
$\domwts_X$ intersects the interior of $\domwts_{G,P}$,
since otherwise $[S,S]$ would be larger.
Applying Lemma~\ref{lemortho},
we conclude 
$$
\cowts_{M^0_S/[M_S,M_S]}\cap\poscowts_{G,P}=\langle0\rangle.
$$

Conversely, suppose $S\subset G$ is a horospherical subgroup
with normalizer $P\subset G$.
We define $X$ to be the spectrum of the ring $\BC[X]$ of
(right) $S$-invariants in the
ring of regular functions $\BC[G]$. Then $\BC[X]$ is finitely-generated,
since $S$ contains the unipotent
radical of a Borel subgroup of $G$.
We have 
$\domwts_X\subset \domwts_{G,P}$, since otherwise $[S,S]$ would be smaller.
Suppose 
$$
\cowts_{M^0_S/[M_S,M_S]}\cap\poscowts_{G,P}=\langle0\rangle.
$$
Applying Lemma~\ref{lemortho}, we conclude that
 $\domwts_X$ intersects the interior of $\domwts_{G,P}$.
Therefore $S/[S,S]$ consists of exactly those elements of $P/[P,P]$ annhilated by $\domwts_X$,
and so $S$ is the generic stabilizer of $X$.
\end{proof}

\end{subsection}


\begin{subsection}{Canonical affine closure}
Let $S\subset G$ be the generic stabilizer of an affine horospherical
$G$-variety $X$ containing a dense $G$-orbit.
Let $\BC[G]$ be the ring of regular functions on $G$,
and let $\BC[G]^S\subset\BC[G]$ be the (right) $S$-invariants.
We call the affine variety
$$
\ol{G/S}=\Spec(\BC[G]^S)
$$
the canonical affine closure of $G/U$.
We have the natural map
$$
\overline{G/S}\to X
$$
corresponding to the 
restriction map 
$$
\BC[X]\to\BC[G/S]\simeq\BC[G]^S.
$$
Since $S$ is horospherical,
the ring $\BC[G]^S$ is simple and graded, and so the affine variety $\ol{G/S}$
is spherical and horospherical.

Although we do not use the following, it clarifies the relation
between $X$ and the canonical affine closure $\ol{G/S}$.

\begin{prop} Let $X$ be an affine horospherical $G$-variety
containing a dense $G$-orbit and generic stabilizer $S\subset G$.
The semigroup 
$
\domwts_{\ol{G/S}}\subset\wts_G
$ 
is the intersection of the dominant weights $\domwts_G\subset\wts_G$ with
the group generated by the semigroup $\domwts_X\subset\wts_G$. 
\end{prop}

\begin{proof}
Let $P\subset G$ be the normalizer of $S\subset G$.
The intersection of $\domwts_G$ and the group generated by $\domwts_X$
consists of exactly those weights in $\domwts_{G,P}$
that annhilate $S/[S,S]$.
\end{proof}

\end{subsection}

\end{section}


\begin{section}{Ind-stacks}

As usual, let $C$ be a smooth complete complex algebraic curve.

\begin{subsection}{Labellings}
Fix a pair $(\lat,\poslat)$ of a lattice $\lat$ and a semigroup $\poslat\subset\lat$.
We shall apply the following to the pair
$(\cowts_{M/[M,M]},\poscowts_{G,P})$.

For $\theta^\pos\in\poslat$, we write $\fU(\theta^\pos)$ 
for a decomposition 
$$
\theta^\pos=\sum_m n_m\theta^\pos_m
$$
where $\theta^\pos_m\in\poslat\setminus \{0\}$ 
are pairwise distinct and $n_m$ are positive integers.

For $\theta^\pos\in\poslat$, and
a decomposition $\fU(\theta^\pos)$,
we write $C^{\fU(\theta^\pos)}$ for 
the partially symmetrized
power $\prod_m C^{(n_m)}$ 
of the curve $C$. 
We write $C^{\fU(\theta^\pos)}_0\subset C^{\fU(\theta^\pos)}$ for the complement
of the diagonal divisor.

For $\Theta$ a pair $(\theta,\fU(\theta^\pos))$ consisting of $\theta\in\lat$,
and $\fU(\theta^\pos)$ a decomposition of $\theta^\pos\in\poslat$,
we write $C^\Theta$ for the product $C\times C^{\fU(\theta^\pos)}$.
We write $C_0^{\Theta}\subset C^{\Theta}$ for the complement
of the diagonal divisor.
Although $C^\Theta$ is independent of $\theta$, it is notationally convenient 
to denote it as we do.

\end{subsection}


\begin{subsection}{$\overline{\on{\mathbf{ Ind-stack}}}$ associated to parabolic subgroup}

Fix a parabolic subgroup $P\subset G$, and 
let $M$ be its Levi quotient $P/U(P)$.
For our application, $P$ will be the 
normalizer of the generic stabilizer $S\subset G$ 
of an irreducible affine horospherical $G$-variety.

Let $\bBunP$ be the ind-stack
that classifies data
$$
(c\in C,\CP_G\in \Bun_G,\CP_{M/[M,M]}\in\Bun_{M/[M,M]},\sigma:\CP_{M/[M,M]}|_{C\setminus c}\to 
\CP_G{\gtimes}\overline{G/[P,P]}|_{C\setminus c})
$$
where $\sigma$ is an $M/[M,M]$-equivariant section which factors
$$
\sigma|_{C'}:\CP_{M/[M,M]}|_{C'}\to 
\CP_G{\gtimes}{G/[P,P]}|_{C'}\to
\CP_G{\gtimes}\overline{G/[P,P]}|_{C'}
$$
for some open curve $C'\subset C\setminus c$.

\begin{subsubsection}{Stratification}
Let $\Theta$ be a pair $(\theta,\fU(\theta^\pos))$, with ${\theta}\in\cowts_{M/[M,M]}$,
and $\theta^\pos\in\poscowts_{G,P}$.
We recall that we have a locally closed embedding
$$
j_{\Theta}:\Bun_P\times C^\Theta_0\to\bBunP
$$ 
defined by
$$
j_\Theta(\CP_P, (c,\sum_{m,n} \theta^\pos_m\cdot c_{m,n}))=
(c,\CP_P\overset{P}{\times} G,\CP_{P}\overset{P}{\times}[P,P]
(-\theta\cdot c-\sum_{m,n} \theta^\pos_m\cdot c_{m,n}),\sigma)
$$ 
where $\sigma$ is the natural map
$$
\CP_{P}\overset{P}{\times}[P,P](-\theta\cdot c-\sum_{m,n} \theta^\pos_m\cdot c_{m,n})
|_{C\setminus c}\to 
\CP_P\overset{P}{\times} G{\gtimes}\overline{G/[P,P]}|_{C\setminus c}
$$
induced by the inclusion
$$
\CP_{P}\overset{P}{\times}\ol{P/[P,P]}
\subset  
\CP_P\overset{P}{\times} \overline{G/[P,P]}
\simeq
\CP_P\overset{P}{\times} G{\gtimes}\overline{G/[P,P]}.
$$
The following is an ind-version of 
\cite[Propositions 6.1.2 \& 6.1.3]{BG02}, or \cite[Proposition 1.5]{BFGMsel02},
and we leave the proof to the reader.

\begin{prop}\label{propbBunPstrat}
Let $\Theta$ be a pair $(\theta,\fU(\theta^\pos))$, with ${\theta}\in\cowts_{M/[M,M]}$,
and $\theta^\pos\in\poscowts_{G,P}$.

Every closed point of $\bBunP$ belongs to the image of a unique $j_\Theta$.
\end{prop}

For $\Theta$ a pair $(\theta,\fU(\theta^\pos))$, with ${\theta}\in\cowts_{M/[M,M]}$,
and $\theta^\pos\in\poscowts_{G,P}$,
we write $\bBunP^\Theta\subset \bBunP$  for the image of $j_\Theta$,
and $\bBunP^{\leq\Theta}\subset\bBunP$ for the closure
of $\bBunP^\Theta\subset\bBunP$.

For $\Theta$ a pair $(\theta,\fU(0))$, with ${\theta}\in\cowts_{M/[M,M]}$, 
the substack $\bBunP^\Theta\subset\bBunP$ 
classifies data
$(c,\CP_G,\CP_{M/[M,M]},\sigma)$
for which the map 
$$
\CP_{M/[M,M]}(\theta\cdot c)|_{C\setminus c}
\stackrel{\sigma}{\to} 
\CP_G{\gtimes}\overline{G/[P,P]}|_{C\setminus c}
$$
extends to a holomorphic map 
$$
\CP_{M/[M,M]}(\theta\cdot c)
\stackrel{\sigma}{\to} 
\CP_G{\gtimes}\overline{G/[P,P]}
$$
which factors
$$
\CP_{M/[M,M]}(\theta\cdot c)
\stackrel{\sigma}{\to} 
\CP_G{\gtimes}{G/[P,P]}\to
\CP_G{\gtimes}\overline{G/[P,P]}.
$$
In this case,
we write $j_\theta$ in place of $j_\Theta$,
$\bBunP^\theta$ in place of $\bBunP^\Theta$,
and $\bBunP^{\leq\theta}$ in place of $\bBunP^{\leq\Theta}$.
For example,
$\bBunP^{\leq 0}\subset\bBunP$ is the closure of the canonical embedding 
$$j_0:\Bun_P\times C\to\bBunP.
$$

\end{subsubsection}

\end{subsection}


\begin{subsection}{$\widetilde{\on{\bf Ind-stack}}$ associated to parabolic subgroup}
Fix a parabolic subgroup $P\subset G$, and let $M$ be its Levi quotient $P/U(P)$.
As usual,
for our application, $P$ will be the 
normalizer of the generic stabilizer $S\subset G$ 
of an irreducible affine horospherical $G$-variety.

Let $\tBunP$ be the ind-stack
that classifies data
$$
(c\in C,\CP_G\in\Bun_G,\CP_{M}\in\Bun_M,\sigma:\CP_{M}|_{C\setminus c}\to 
\CP_G{\gtimes}\overline{G/U(P)}|_{C\setminus c})
$$
where 
$\sigma$ is an $M$-equivariant section
which factors
$$
\sigma|_{C'}:\CP_{M}|_{C'}\to 
\CP_G{\gtimes}{G/U(P)}|_{C'}\to
\CP_G{\gtimes}\overline{G/U(P)}|_{C'}
$$
for some open curve $C'\subset C\setminus c$.

\begin{subsubsection}{Stratification}

For 
$\theta^\pos\in\poscowts_{G,P}$, 
we write $\tilde\fU(\theta^\pos)$ 
for a collection of (not necessarily distinct) elements $\tilde\theta^\pos_{m}\in\tilde\cowts^\pos_{G,P}\setminus\{0\}$ 
such that
$$
\theta^\pos=\sum_{m} r(\tilde\theta^\pos_{m}).
$$
We write $r(\tilde\fU(\theta^\pos))$ for the decomposition such a collection defines.

Let $\tilde\Theta$ be a pair $(\tilde\theta,\tilde\fU(\theta^\pos))$ with 
$\tilde\theta\in\domcowts_M$, and $\theta^\pos\in\cowts^\pos_{G,P}$,
and let $\Theta$ be the associated pair $(r(\tilde\theta),r(\tilde\fU(\theta^\pos)))$.
We define the Hecke ind-stack
$$
\H^{\tilde\Theta}_{M,0}
\to
C_0^\Theta
$$
to be that with fiber over $(c,c_{\fU(\theta^\pos)})\in C_0^\Theta$, where $c_{\fU(\theta^\pos)}=\sum_{m} r(\tilde \theta^\pos_m)\cdot c_{m}$,
the fiber product
$$
\H_{M}^{\tilde\theta}|_c
\timesBunM\prod_{{\Bun_M}} {\H}_{M}^{\tilde\theta^\pos_m}|_{c_{m}}.
$$

The following is an ind-version of \cite[Proposition 6.2.5]{BG02}, or 
\cite[Proposition 1.9]{BFGMsel02},
and we leave the proof to the reader.

\begin{prop}\label{pttobfiber}
Let $\tilde\Theta$ be a pair $(\tilde\theta,\tilde\fU(\theta^\pos))$ with 
$\tilde\theta\in\domcowts_M$, and $\theta^\pos\in\cowts^\pos_{G,P}$.

On the level of reduced ind-stacks,
there is a locally closed embedding 
$$
j_{\tilde\Theta}:\Bun_P\underset{\Bun_M}{\times} {\H_{M,0}^{\tilde\Theta}}\to\tBunP.
$$
Every closed point of $\tBunP$ belongs to the image of a unique $j_{\tilde\Theta}$.

\end{prop}

For $\tilde\Theta$ a pair $(\tilde\theta,\tilde\fU(\theta^\pos))$, with 
${\tilde\theta}\in\domcowts_{M}$,
and $\theta^\pos\in\poscowts_{G,P}$,
we write $\tBunP^{\tilde\Theta}\subset \tBunP$  for the image of $j_{\tilde\Theta}$,
and $\tBunP^{\leq\tilde\Theta}\subset\tBunP$ for the closure
of $\tBunP^{\tilde\Theta}\subset\tBunP$.

For $\tilde\Theta$ a pair $(\tilde\theta,\tilde\fU(0))$, with 
$\tilde\theta\in\domcowts_M$, we write 
$j_{\tilde\theta}$ in place of $j_{\tilde\Theta}$,
$\tBunP^{\tilde\theta}$ in place of $\tBunP^{\tilde\Theta}$,
and $\tBunP^{\leq\tilde\theta}$ in place of $\tBunP^{\leq\tilde\Theta}$
For example, 
$\tBunP^{\leq 0}$ is the closure of the canonical embedding
$$
j_{\tilde 0}:\Bun_P\times C\to\tBunP.
$$ 

\end{subsubsection}

\end{subsection}


\begin{subsection}{$\overline{\on{\bf Ind-stack}}$ associated to generic stabilizer}
Let $X$ be an irreducible affine horospherical $G$-variety with generic stabilizer $S\subset G$.
Recall 
that the normalizer of $S$ is a parabolic subgroup $P\subset G$
with the same derived group $[P,P]=[S,S]$ and unipotent radical $U(P)=U(S)$.
Let $M$ be the Levi quotient $P/U(P)$, and let $M_S$  be the Levi quotient
$S/U(S)$.

Let $\bZcan$ be the ind-stack that classifies data
$$
(c\in C,\CP_G\in\Bun_G,\CP_{M_S/[M_S,M_S]}\in\Bun_{M_S/[M_S,M_S]},
\sigma:\CP_{M_S/[M_S,M_S]}|_{C\setminus c}\to 
\CP_G{\gtimes}\overline{G/[S,S]}|_{C\setminus c})
$$
where $\sigma$ is an $M_S/[M_S,M_S]$-equivariant section
which factors
$$
\sigma|_{C'}:\CP_{M_S/[M_S,M_S]}|_{C'}\to 
\CP_G{\gtimes}{G/[S,S]}|_{C'}\to
\CP_G{\gtimes}\overline{G/[S,S]}|_{C'}
$$
for some open curve $C'\subset C\setminus c$.

The following is immediate from the definitions.

\begin{prop}\label{propbZcanCart}
The diagram
$$
\begin{array}{ccc}
\bZcan & \to & \bBunP \\
\downarrow & & \downarrow \\
\Bun_{M_S/[M_S,M_S]} & \to & \Bun_{M/[M,M]} \\
\end{array}
$$
is Cartesian.
\end{prop}

\begin{subsubsection}{Stratification}
Let $\Theta$ be a pair $(\theta,\fU(\theta^\pos))$, with $\theta\in\cowts_{M/[M,M]}$,
and $\theta^\pos\in\poscowts_{G,P}$.

We write $\bZcan^\Theta\subset \bZcan$ for the substack
which completes the Cartesian diagram
$$
\begin{array}{ccc}
\bZcan^\Theta & \to & \bBunP^\Theta \\
\downarrow & & \downarrow \\
\Bun_{M_S/[M_S,M_S]} & \to & \Bun_{M/[M,M]}, \\
\end{array}
$$
and $\bZcan^{\leq\Theta}\subset \bZcan$ for the closure of
$\bZcan^\Theta\subset \bZcan$.

For $\Theta$ a pair $(\theta,\fU(0))$, with ${\theta}\in\cowts_{M/[M,M]}$, 
we write
$\bZcan^{\theta}$ in place of $\bZcan^\Theta$, and
$\bZcan^{\leq\theta}$ in place of 
$\bZcan^{\leq\Theta}$.
For example,
$\bZcan^{\leq 0}$ is the closure of the canonical embedding 
$$
\Bun_S\times C\subset\bZcan.
$$

\end{subsubsection}

\end{subsection}


\begin{subsection}{Naive ind-stack associated to $X$}\label{sectnaive}
Let $X$ be an affine horospherical $G$-variety with dense $G$-orbit $\openX\subset X$
and generic stabilizer $S\subset G$.

Let $Z$ be the ind-stack that classifies data
$$
(c\in C, \CP_G\in\Bun_G,\sigma:{C\setminus c}\to 
\CP_G{\gtimes}X|_{C\setminus c})
$$
where $\sigma$ is a section which
factors
$$
\sigma|_{C'}:C'\to 
\CP_G{\gtimes}\openX|_{C'}\to
\CP_G{\gtimes}{X}|_{C'}
$$
for some open curve $C'\subset C\setminus c$.

For the canonical affine closure $\overline{G/S}$,
we write $\Z_\can$ for the corresponding ind-stack.

We call the ind-stack $\Z$ naive, since there is no auxilliary bundle
in its definition:
it classifies honest sections.
Let $\stZ$ be the ind-stack 
that classifies data
$$
(c\in C,\CP_G\in\Bun_G,\CP_{M/M_S}\in\Bun_{M/M_S},
\sigma: \CP_{M/M_S}|_{C\setminus c}\to 
\CP_G{\gtimes}X|_{C\setminus c})
$$
where $\sigma$ is an $M/M_S$-equivariant section
which factors
$$
\sigma|_{C'}:\CP_{M/M_S}|_{C'}\to 
\CP_G{\gtimes}{\openX}|_{C'}\to
\CP_G{\gtimes}{X}|_{C'}
$$
for some open curve $C'\subset C\setminus c$.
Here as usual, we write $M$ for the Levi quotient $P/U(P)$ of the normalizer $P\subset G$
of the generic stabilizer $S\subset G$,
and $M_S$  for the Levi quotient
$S/U(S)$.

For the canonical affine closure $\overline{G/S}$, 
we write $\stZcan$ for the corresponding ind-stack.

The following analogue of Proposition~\ref{propbZcanCart}
is immediate from the definitions.

\begin{prop}
The diagram
$$
\begin{array}{ccc}
\Z & \to & \stZ \\
\downarrow & & \downarrow \\
\Bun_{\langle1\rangle} & \to & \Bun_{M/M_S} \\
\end{array}
$$
is Cartesian.
\end{prop}

%
%
\begin{subsubsection}{Stratification}
We shall content ourselves here with defining the
substacks of the naive ind-stack $\Z$ which appear in our main theorem.
(See \cite{GNcombo06} for a different perspective involving a completely local definition.)
Recall that we write $A$ for the quotient torus $P/S$,
and $\cowts_A$ for its coweight lattice. 
Similarly, for the identity component $S^0\subset S$,
we write $A_0$ for the quotient torus $P/S^0$,
and $\cowts_{A_0}$ for its coweight lattice. 
The natural map $A_0\to A$ provides
an inclusion of coweight lattices $\cowts_{A_0}\to \cowts_A$.
For $\kappa\in\cowts_{A}$,
we shall define a closed substack $\Z^{\leq\kappa}\subset\Z$.
When $\kappa\in\cowts_{A_0}$,
the closed substack $\Z^{\leq\kappa}\subset\Z$
appears in our main theorem.

For $\kappa\in\cowts_A$, 
let
$\stZ^{\kappa}\subset \stZ$ be the locally closed substack 
that classifies data
$(c,\CP_G,\CP_{M/M_S},\sigma)$
for which the natural map 
$$
\CP_{M/M_S}(\kappa\cdot c)|_{C\setminus c}
\stackrel{\sigma}{\to} 
\CP_G{\gtimes} X|_{C\setminus c}
$$
extends to a holomorphic map 
$$
\CP_{M/M_S}(\kappa\cdot c)
\stackrel{\sigma}{\to} 
\CP_G{\gtimes} X
$$
which factors
$$
\CP_{M/M_S}(\kappa\cdot c)
\stackrel{\sigma}{\to} 
\CP_G{\gtimes}{\openX}\to
\CP_G{\gtimes}X.
$$
We write $\stZ^{\leq\kappa}\subset\stZ$ for the closure of 
$\stZ^{\kappa}\subset\stZ$.

For $\kappa\in\cowts_A$, 
let 
$\Z^{\kappa}\subset\Z$ be the locally closed substack completing the Cartesian diagram
$$
\begin{array}{ccc}
\Z^{\kappa} & \to & \stZ^{\kappa} \\
\downarrow & & \downarrow \\
\Bun_{\langle1\rangle} & \to & \Bun_{M/M_S}. \\
\end{array}
$$
We write $\Z^{\leq\kappa}\subset\Z$ for the closure of 
$\Z^{\kappa}\subset\Z$.

\end{subsubsection}

\end{subsection}

\end{section}


\begin{section}{Maps}

\begin{subsection}{The map $\fr:\tBunP\to\bBunP$}

Let $\Theta$ be a pair $(\theta,\fU(\theta^\pos))$, with ${\theta}\in\cowts_{M/[M,M]}$,
and 
$\theta^\pos\in\poscowts_{G,P}$.
and $\fU(\theta^\pos)$ a decomposition $\theta^\pos=\sum_{m} n_m\theta^\pos_m$.
Let $\tBunP^\Theta\subset\tBunP$ be the inverse image of $\bBunP^\Theta\subset\bBunP$ 
under the natural map 
$$
\fr:\tBunP\to\bBunP.
$$ 
We would like to describe the fibers
of the restriction of $\fr$ to the substack $\tBunP^\Theta\subset\tBunP$.

First, we define the Hecke ind-substack
$$
{\H}_{M}^{\flat(\theta)}\subset
{\H}_{M}
$$ 
to be the union of the spherical Hecke substacks
$$
{\H}_{M}^{\mu}\subset
{\H}_{M}, 
$$
for $\mu\in\domcowts_{M}$ such that $r(\mu)=\theta$.

Second,
if there exists $\tilde\mu^\pos\in\tilde{\Lambda}^\pos_{G,P}$ such that
$r(\tilde\mu^\pos)=\theta^\pos$, we
define the Hecke substack
$$
{\H}_{M}^{\flat(\theta^\pos)}\subset
{\H}_{M}
$$ 
to be the union of the spherical Hecke substacks
$$
{\H}_{M}^{\tilde\mu^\pos}\subset
{\H}_{M},
$$ for $\tilde\mu^\pos\in\tilde{\cowts}^\pos_{G,P}$ such that 
$r(\tilde\mu^\pos)=\theta_m^\pos$.

Finally, we define the Hecke ind-stack
$$
\H_{M,0}^{\flat(\Theta)}
\to
C_0^\Theta
$$
to be that with fiber over $(c,c_{\fU(\theta^\pos)})\in C_0^\Theta$, where $c_{\fU(\theta^\pos)}=\sum_{m,n}  \theta^\pos_m\cdot c_{m,n}$,
the fiber product
$$
\H_{M}^{\flat(\theta)}|_c
\timesBunM\prod_{{\Bun_M}} {\H}_{M}^{\flat(\theta^\pos_m)}|_{c_{m,n}}.
$$

The following is an ind-version of \cite[Proposition 6.2.5]{BG02}, or 
\cite[Proposition 1.9]{BFGMsel02},
and we leave the proof to the reader.
It is also immediately implied by Proposition~\ref{pttobfiber}.

\begin{prop}\label{proptBunPstrat}
Let $\Theta$ be a pair $(\theta,\fU(\theta^\pos))$, with $\theta\in\cowts_{M/[M,M]}$,
$\theta^\pos\in\poscowts_{G,P}$,
and $\fU(\theta^\pos)$ a decomposition $\theta^\pos=\sum_{m} n_m\theta^\pos_m$.

If for all $m$ there exists $\tilde\mu_m^\pos\in\tilde\cowts^\pos_{G,P}$ such that
$r(\tilde\mu_m^\pos)=\theta^\pos_m$,
then on the level of reduced stacks there is a canonical isomorphism
$$
\tBunP^\Theta\simeq\Bun_P\underset{\Bun_M}{\times}\H_{M,0}^{\flat(\Theta)}
$$
such that the following diagram commutes
$$
\begin{array}{ccc}
\tBunP^\Theta & \simeq & \Bun_P\underset{\Bun_M}{\times}\H_{M,0}^{\flat(\Theta)} \\
\downarrow & & \downarrow \\
\bBunP^\Theta & \simeq & \Bun_P\times C_0^\Theta \\
\end{array}
$$
where the right hand side is the obvious projection.

If there is an $m$ such that $\theta^\pos_m$ is not equal to $r(\tilde\mu^\pos)$,
for any $\tilde\mu^\pos\in\tilde\cowts^\pos_{G,P}$,
then 
$\tBunP^\Theta$ is empty.
\end{prop}

\end{subsection}


\begin{subsection}{The map $\fp:\bZcan\to\Zcan$.}

Let $X$ be an irreducible affine horospherical $G$-variety
with generic stabilizer $S\subset G$.
Recall that the normalizer of a horospherical subgroup $S\subset G$
is a parabolic subgroup $P\subset G$ with the same derived group $[P,P]=[S,S]$
and unipotent radical $U(P)=U(S)$.
We write $M$ for the Levi quotient $P/U(P)$, $M_S$ for the Levi quotient
$S/U(S)$, and $M^0_S$ for the identity component of $M_S$.
We write $A$ for the quotient torus $P/S$,
and $\cowts_A$ for its coweight lattice. 
Similarly, for the identity component $S^0\subset S$,
we write $A_0$ for the quotient torus $P/S^0$,
and $\cowts_{A_0}$ for its coweight lattice. 
The natural map $M/[M,M]\to A_0$ induces a surjection of coweight lattices
$\cowts_{M/[M,M]}\to\cowts_{A_0}$ which we denote by $p$.
The kernel of $p$ is the coweight lattice $\cowts_{M^0_S/[M_S,M_S]}$.
(Note that the component group of $M_S$ is abelian.)

Associated to the canonical affine closure $\overline{G/S}$,
we have a Cartesian diagram of ind-stacks
$$
\begin{array}{ccc}
\bZcan & \to & \bBunP \\
\fp\downarrow & & \downarrow\fp \\
\Zcan & \to & \stZcan \\
\end{array}
$$
We would like to describe some properties of the vertical maps.

\begin{prop}
The map $\fp:\bBunP\to\stZcan$ is ind-finite.

For $\theta\in\cowts_{M/[M,M]}$,
its restriction to $\bBunP^\theta$ is an embedding with
image $\stZ_\can^{p(\theta)}$, and
its restriction to $\bBunP^{\leq\theta}$ is finite with
image $\stZ_\can^{\leq p(\theta)}$.
\end{prop}

\begin{proof}
For a point $(c,\CP_G,\CP_{M/[M,M]},\overline\sigma)\in\bBunP$,
we write $(c,\CP_G,\CP_{M/M_S},\sigma)\in\stZcan$ for its
image under $\fp$.
Observe that for $\theta\in\cowts_{M/[M,M]}$,
the point $(c,\CP_G,\CP_{M/[M,M]}(\theta\cdot c),\overline\sigma)\in\bBunP$
maps to
$(c,\CP_G,\CP_{M/M_S}(p(\theta)\cdot c),\sigma)\in\stZcan$ under $\fp$.
Therefore to prove the proposition, it suffices to show that
the restriction of $\fp$ to the canonical embedding
$\Bun_P\subset\bBunP$ is an embedding
with image the canonical embedding $\Bun_P\subset\stZcan$,
and its restriction to $\bBunP^{\leq 0}$ is a finite map with
image $\stZ_{\can}^{\leq 0}$.
The first assertion is immediate from the definitions.
To prove the second,
recall that by \cite[Proposition 1.3.6]{BG02},
$\bBunP$ is proper over $\Bun_G$, and so the map $\fp$ is proper
since it respects the projection to $\Bun_G$.
Therefore 
it suffices to check that the fibers 
over closed points 
of the restriction of $\fp$ to $\bBunP^{\leq 0}$
are finite.

Let $\Theta$ be a pair $(0,\fU(\theta^\pos))$, with 
$\theta^\pos\in\poscowts_{G,P}$.
The stack
$\bBunP^\Theta$ classifies data
$$
(c,\CP_P,c_\Theta,\CP_{M/[M,M]})
$$
together with an isomorphism
$$
\alpha:\CP_P\overset{P}{\times}P/[P,P]\simeq\CP_{M/[M,M]}(c_\Theta).
$$
The fiber of $\fp$ through such a point classifies data
$$
(\CP_P,c'_{\Theta'},\CP'_{M/[M,M]})
$$
together with an isomorphism
$$
\alpha':\CP_P\overset{P}{\times}P/[P,P]\simeq\CP'_{M/[M,M]}(c'_{\Theta'})
$$
such that the labelling
$c_\Phi=c_{\Theta}-c'_{\Theta'}$ takes values in $\cowts_{M^0_S/[M_S,M_S]}$.
Therefore 
we need only check that for
$\theta^\pos\in\poscowts_{G,P}$, there are only a finite number
of $\phi\in\cowts_{M^0_S/[M_S,M_S]}$ such that
$\theta^\pos+\phi\in\poscowts_{G,P}$.
By Proposition~\ref{propstab}, the lattice $\cowts_{M^0_S/[M_S,M_S]}$ intersects
the semigroup $\poscowts_{G,P}$ only at $0$.
Since $\poscowts_{G,P}$ is finitely-generated, this implies
that for $\theta^\pos\in\cowts_{M/[M,M]}$, the coset
$\theta^\pos+\cowts_{M^0_S/[M_S,M_S]}$ intersects $\poscowts_{G,P}$
in a finite set. 
\end{proof}

\begin{cor}\label{corsmall}
The map $\fp:\bZcan\to\Zcan$ is ind-finite.

For $\theta\in\cowts_{M/[M,M]}$,
its restriction to $\bZcan^\theta$ is an embedding with
image $\Z_\can^{p(\theta)}$, and
its restriction to $\bZcan^{\leq\theta}$ is finite with
image $\Z_\can^{\leq p(\theta)}$.

\end{cor}

\end{subsection}


\begin{subsection}{The map $\fs:\Zcan\to \Z$.}
Let $X$ be an affine horospherical variety
with dense $G$-orbit $\openX\subset X$ and generic stabilizer $S\subset G$.

Associated to the natural map $\ol{G/S}\to X$,
we have a Cartesian diagram of ind-stacks
$$
\begin{array}{ccc}
\Zcan & \to & \stZcan \\
\fs\downarrow & & \downarrow \fs\\
\Z & \to & \stZ. \\
\end{array}
$$
We would like to describe some properties of the vertical maps.

\begin{prop}\label{propsembed}
The map $\fs:\stZcan\to\stZ$ is a closed embedding.

For $\kappa\in\cowts_{A}$,
its restriction to $\stZ_{\can}^{\kappa}$ is an embedding with
image $\stZ^{\kappa}$, and
its restriction to $\stZ_{\can}^{\leq\kappa}$ is a closed embedding with
image $\stZ^{\leq\kappa}$.
\end{prop}

\begin{proof}
First note that $\fs$ is injective on scheme-valued points since for 
$(c,\CP_G,\CP_{M/M_S}\sigma)\in\stZcan$,
the map 
$$\sigma:\CP_{M/M_S}|_{C\setminus c}\to\CP_G\gtimes \overline{G/S}|_{C\setminus c}$$
factors
$$
\sigma|_{C'}:\CP_{M/M_S}|_{C'}\to\CP_G \gtimes{G/S}|_{C'}\to\CP_G \gtimes\overline{G/S}|_{C'},
$$
for some open curve $C'\subset C\setminus c$,
and the map $\overline{G/S} \to X$ restricted to $G/S $ is an embedding.

Now to see $\fs$ is a closed embedding,
it suffices to check that $\fs$ satisfies the valuative criterion of properness.
Let $D=\Spec\BC[[t]]$ be the disk, and $D^\times=\Spec\BC((t))$ the punctured disk.
Let $f:D\to \Z$ be a map with a partial lift $F^\times:D^\times\to\Zcan$.
Let $\CP_G^f$ be the $D$-family of $G$-bundles defined by $f$,
and let $\CP_{M/M_S}^f$ be the $D$-family of $M/M_S$-bundles defined by $f$.
We must check that any partial lift 
$$
\Sigma^\times:\CP_{M/M_S}^f|_{(C\setminus c)\times D^\times}
\to\CP_G^f \gtimes \overline{G/S}|_{(C\setminus c)\times D^\times}
$$
of a map 
$$
\sigma:\CP_{M/M_S}^f|_{(C\setminus c)\times D}\to\CP_G^f \gtimes X|_{(C\setminus c)\times D}
$$
which factors
$$
\sigma|_{C'\times D}:\CP_{M/M_S}^f|_{C'\times D}\to\CP_G^f \gtimes{G/S}|_{C'\times D}
\to\CP_G^f \gtimes X|_{C'\times D},
$$
for some open curve $C'\subset C\setminus c$,
extends to $(C\setminus c)\times D$.
Since $\overline{G/S}\to X$ restricted to $G/S$ 
is an embedding with image $G/S$, 
we may lift
$\sigma|_{C'\times D}$ 
to extend $\Sigma^\times$ to $C'\times D$.
But then $\Sigma^\times$ extends completely since $\CP_{M/M_S}^f|_{(C\setminus c)\times D}$ is normal and
the complement of  $\CP_{M/M_S}^f|_{C'\times D}$ is of codimension $2$.

Finally, 
for a point $(c,\CP_G,\CP_{M/M_S},\sigma_\can)\in\stZcan $,
we write $(c,\CP_G,\CP_{M/M_S},\sigma)\in\stZ$ for its
image under $\fs$.
Observe that for $\kappa\in\cowts_{A}$, 
the point $(c,\CP_G,\CP_{M/M_S}(\kappa\cdot c),\sigma_\can)\in\stZcan $
maps to
$(c,\CP_G,\CP_{M/M_S}(\kappa\cdot c),\sigma)\in\stZ$ under $\fs$.
Therefore to complete the proof of the proposition,
it suffices to show that
the restriction of $\fs$ to the canonical embedding
$\Bun_S\times C \subset\stZcan $ has image the canonical embedding 
$\Bun_S\times C \subset\stZ$.
This is immediate from the definitions.
\end{proof}

\begin{cor}\label{corsembed}
The map $\fs:\Zcan\to\Z$ is a closed embedding.

For $\kappa\in\cowts_{A}$,
its restriction to $\Z_{\can}^{\kappa}$ is an embedding with
image $\Z^{\kappa}$, and
its restriction to $\Z_{\can}^{\leq\kappa}$ is a closed embedding with
image $\Z^{\leq\kappa}$.
\end{cor}
\end{subsection}

\end{section}


\begin{section}{Convolution}\label{sconv}
Let $X$ be an affine horospherical $G$-variety with dense $G$-orbit $\openX\subset X$
and generic stabilizer $S\subset G$.

The following diagram summarizes the ind-stacks and maps under consideration
$$
\begin{array}{ccccccc}
\tBunP & \stackrel{\fr}{\to} & \bBunP  & \stackrel{\fp}{\to} & \stZcan   && \\
&& \uparrow\fk && \uparrow\fk && \\
&& \bZcan  & \stackrel{\fp}{\to} & \Zcan  & \stackrel{\fs}{\to} & \Z.\\
\end{array}
$$
Each of the ind-stacks of the diagram projects to $C\times\Bun_G$, and 
the maps of the diagram commute with the projections.

Let $\CZ$ be any one of the ind-stacks from the diagram, 
and form the diagram
$$
\begin{array}{ccccc}
\CZ & \stackrel{{}\hlG}{\leftarrow} & \H_G\underset{Bun_G\times C}{\times}\CZ
& \stackrel{{}\hrG}{\to} & \CZ \\
\downarrow & & \downarrow & & \downarrow \\
\Bun_G & \stackrel{{}\hlG}{\leftarrow} & \H_G & \stackrel{{}\hrG}{\to} & \Bun_G
\end{array}
$$
in which each square is Cartesian.

For $\lambda\in\domcowts_G$, 
we define the convolution functor 
$$
H^\lambda_G :\Sh(\CZ)\to\Sh(\CZ)
$$ 
on an object $\CF\in\Sh(\CZ)$ to be
$$
H^\lambda_G (\CF)={}{\hlG}_!(\CA^\lambda_G\tboxtimes\CF)^r
$$ 
where $(\CA^\lambda_G\tboxtimes\CF)^r$
is the twisted product defined with respect to ${}\hrG$, and $\CA_G^\lambda$
is the simple spherical sheaf on the fibers of $\hrG$ corresponding to $\lambda$.
(See Section~\ref{sbandh} for more on the twisted product and spherical sheaf.)


\begin{subsection}{Convolution on $\tBunP$}

Recall that for a reductive group $H$, and $\lambda\in\domcowts_H$, 
we write $V_\chH^\lambda$ for the irreducible
representation of the dual group $\chH$ of highest weight $\lambda$. 

We shall deduce our results from the following.

\begin{thm}\label{thmtBunPconv}\cite[Theorem 4.1.5]{BG02}.
For $\lambda\in\domcowts_G$, there is a canonical isomorphism
$$
H^\lambda_G (\IC_{\tBunP}^{\leq 0})\simeq
\sum_{\mu\in\domcowts_M}
\IC_{\tBunP}^{\leq\mu} \otimes \Hom_{\check M} (V_\chM^\mu,V_\chG^\lambda).
$$
\end{thm}

\end{subsection}


\begin{subsection}{Convolution on $\bBunP$}
Recall that $r:\cowts_M\to\cowts_{M/[M,M]}$ denotes the natural projection,
$2\chrho_M$ the sum of the positive roots of $M$,
and $\langle2\chrho_M,\mu\rangle$ the natural pairing, for $\mu\in\cowts_M$.

\begin{thm}\label{thmbBunPconv}
For $\lambda\in\domcowts_G$, there is an isomorphism
$$
H^\lambda_G (\IC_{\bBunP}^{\leq 0})\simeq
\sum_{\theta\in\cowts_{M/[M,M]}}
\sum_{
{\mu\in\cowts_M},
{r(\mu)=\theta} 
}
\IC_{\bBunP}^{\leq\theta} \otimes \Hom_{\chT} (V_\chT^\mu,V_\chG^\lambda)
[\langle2\chrho_M,\mu\rangle].
$$
\end{thm}

\begin{proof}{\em Step 1.}
For the projection
$$
\fr:\tBunP\to\bBunP,
$$
we clearly have
\begin{equation}\label{eqnstep1}
H^\lambda_G (\fr_!\IC_{\tBunP}^{\leq 0})\simeq
\fr_!H^\lambda_G (\IC_{\tBunP}^{\leq 0}).
\end{equation}

Let us first analyze the left hand side of equation~\ref{eqnstep1}.
We may write the pushforward
$\fr_!\IC_{\tBunP}^{\leq 0}$ in the form
$$
\fr_!\IC_{\tBunP}^{\leq 0}\simeq \IC_{\bBunP}^{\leq 0}
\oplus
\CI^{\leq 0}
$$
where $\CI^{\leq 0}\in\Sh(\bBunP)$ is isomorphic to a direct sum
of shifts of sheaves of the form 
$$
\IC_{\bBunP}^{\leq\Theta}, 
\mbox{ for pairs 
$\Theta=(0,\fU(\theta^\pos))$, with 
$\theta^\pos\in\poscowts_{G,P}\setminus\{0\}$.}
$$
The asserted form of $\CI^{\leq0}$ follows from the Decomposition Theorem,
the fact that the restrictions of $\IC_{\tBunP}^{\leq 0}$ to the strata
of $\tBunP$ are constant~\cite[Theorem 1.12]{BFGMsel02}, and the structure
of the map $\fr$ described in Proposition~\ref{proptBunPstrat}.

For any $\eta^\pos\in\poscowts_{G,P}\setminus\{0\}$, and decomposition
$\fU(\eta^\pos)$,
we have the finite map 
$$
\tau_{\fU(\eta^\pos)}:C^{\fU(\eta^\pos)}\times\bBunP\to\bBunP
$$ 
defined by
$$
\tau_{\fU(\eta^\pos)}(\sum_{m,n} \eta^\pos_m\cdot c_{m,n},(c, \CP_G,\CP_{M/[M,M]},\sigma))=
(c,\CP_G,\CP_{M/[M,M]}(-\sum_{m,n}\eta_m^\pos\cdot c_{m,n}),\sigma).
$$ 
Note that for $\eta\in\cowts_{M/[M,M]}$, and $\Theta$ the pair $(\eta,\fU(\eta^\pos))$,
the restriction of $\tau_{\fU(\eta^\pos)}$ provides an isomorphism
$$
\tau_{\fU(\eta^\pos)}:(C^{\fU(\eta^\pos)} \times\bBunP^{\eta})_{0}
\risom\bBunP^{\Theta}
$$
where the domain
completes the Cartesian square
$$
\begin{array}{ccc}
(C^{\fU(\eta^\pos)} \times\bBunP^{\eta})_{0}
& \to &
C^{\fU(\eta^\pos)} \times\bBunP^{\eta} \\
\downarrow & & \downarrow \\
(C^{\fU(\eta^\pos)} \times C)_{0}
& \to & 
C^{\fU(\eta^\pos)} \times C
\end{array}
$$
where as usual
$$
(C^{\fU(\eta^\pos)} \times C)_{0}
 \subset 
C^{\fU(\eta^\pos)} \times C
$$
denotes the complement to the diagonal divisor.

We define the strict full triangulated subcategory of irrelevant sheaves
$$
\IrrelSh(\bBunP)\subset\Sh(\bBunP)
$$
to be that generated by sheaves
of the form 
$$
{\tau_{\fU(\eta^\pos)}}_!(\IC^{\fU(\eta^\pos)}_C\boxtimes\CF)
$$ 
where $\eta^\pos$ runs through $\poscowts_{G,P}\setminus\{0\}$,
$\fU(\eta^\pos)$ runs through decompositions of $\eta^\pos$,
$\IC^{\fU(\eta^\pos)}_C$ denotes the intersection cohomology
sheaf of $C^{\fU(\eta^\pos)}$, 
and $\CF$ runs through objects of $\Sh(\bBunP)$.

\begin{lem}
The sheaf $\CI^{\leq 0}$ is irrelevant.
\end{lem}

\begin{proof}
Let $\Theta$
be a pair $(\theta,\fU(\theta^\pos))$, with $\theta\in\cowts_{M/[M,M]}$, and
$\theta^\pos\in\poscowts_{G,P}\setminus\{0\}$.
Then 
we may realize the sheaf $\IC_{\bBunP}^{\leq\Theta}$ as the pushforward
$$
\IC_{\bBunP}^{\leq\Theta}
\simeq
{\tau_{\fU(\theta^\pos)}}_!
(\IC_{C}^{\fU(\theta^\pos)}\boxtimes\IC_{\bBunP}^\theta)
$$
To see this, we use the isomorphism
$$
\tau_{\fU(\theta^\pos)}:(C^{\fU(\theta^\pos)} \times\bBunP^{\theta})_{0}
\risom\bBunP^{\Theta},
$$
and the fact that $\tau_{\fU(\theta^\pos)}$ is finite.
\end{proof}

\begin{lem}
If $\CE$ is an irrelevant sheaf, then $H^\lambda_G (\CE)$
is an irrelevant sheaf.
\end{lem}

\begin{proof}
Clearly we have a canonical isomorphism
$$
H^\lambda_G ({\tau_{\fU(\eta^\pos)}}_!(\IC_C^{\fU(\eta^\pos)}\boxtimes\CF))\simeq
{\tau_{\fU(\eta^\pos)}}_!(\IC^{\fU(\eta^\pos)}_C\boxtimes H^\lambda_G (\CF)).
$$
\end{proof}

By the preceding lemmas, 
we may write the left hand side of equation~\ref{eqnstep1}
in the form
\begin{equation}\label{eqnstep1left}
H^\lambda_G (\fr_!\IC_{\tBunP}^{\leq 0})
\simeq
H^\lambda_G (\IC_{\bBunP}^{\leq 0})\oplus H^\lambda_G (\CI^{\leq 0})
\end{equation}
where $H^\lambda_G (\CI^{\leq 0})$ is an irrelevant sheaf.

Let us next analyze the right hand side of equation~\ref{eqnstep1}.
By Theorem~\ref{thmtBunPconv}, we have
$$
\fr_! H^\lambda_G (\IC_{\tBunP}^{\leq 0})
\simeq
\sum_{\mu\in\domcowts_M}
\fr_!\IC_{\tBunP}^{\leq\mu} \otimes \Hom_{\check M} (V_\chM^\mu,V_\chG^\lambda).
$$

\begin{lem}
For $\mu\in\domcowts_{M}$, we have
$$
\fr_!\IC_{\tBunP}^{\leq \mu}
\simeq 
\sum_{\nu\in\cowts_M} 
(\IC_{\bBunP}^{\leq r(\mu)}
\oplus \CI^{\leq \mu})\otimes
\Hom_\chT(V^\nu_\chT,V^\mu_\chM)[\langle2\chrho_M,\nu\rangle].
$$
where $\CI^{\leq \mu}$ is isomorphic to a direct sum
of shifts of sheaves of the form 
$$\IC_{\bBunP}^{\leq\Theta}, 
\mbox{ for pairs 
$\Theta=(\theta,\fU(\theta^\pos))$.}
$$
\end{lem}

\begin{proof}

We may form the diagram
$$
\begin{array}{ccccc}
\tBunP & \stackrel{{}\hl_M}{\leftarrow} & \H_M\underset{\Bun_M\times C}{\times}\tBunP
& \stackrel{{}\hr_M}{\to} & \tBunP \\
\downarrow & & \downarrow & & \downarrow \\
\Bun_G & \stackrel{{}\hl_M}{\leftarrow} & \H_G & \stackrel{{}\hr_M}{\to} & \Bun_G
\end{array}
$$
in which each square is Cartesian.
We define the convolution functor 
$$
H^\mu_M :\Sh(\tBunP)\to\Sh(\tBunP)
$$ 
on an object $\CF\in\Sh(\tBunP)$ to be
$$
H^\mu_M (\CF)={}{\hl_M}_!(\CA^\mu_M\tboxtimes\CF)^r
$$ 
where $(\CA^\mu_M\tboxtimes\CF)^r$
is the twisted product defined with respect to ${}\hr_M$, and $\CA_M^\mu$
is the simple spherical sheaf on the fibers of $\hr_M$ corresponding to $\mu$.
Theorem 4.1.3 of~\cite{BG02} provides a canonical isomorphism
$$
H^\mu_M(\IC^{\leq 0}_{\tBunP})\simeq\IC^{\leq\mu}_{\tBunP}.
$$
We also have a commutative diagram
$$
\begin{array}{ccc}
\tBunP & \stackrel{\hl_M}{\leftarrow} & \H_M\underset{\Bun_M\times C}{\times}\tBunP \\
\fr\downarrow & & \downarrow\fr' \\
\bBunP & \stackrel{\hl_{M/[M,M]}}{\leftarrow} & 
\H_{M/[M,M]} \underset{\Bun_{M/[M,M]}\times C}{\times}\bBunP
\end{array}
$$
where the modification map ${\hl_{M/[M,M]}}$ is given by
$$
{\hl_{M/[M,M]}}(\theta,(c, \CP_G,\CP_{M/[M,M]},\sigma))=
(c,\CP_G,\CP_{M/[M,M]}(-\theta\cdot c),\sigma).
$$ 
We conclude that there is an isomorphism
$$
\fr_!\IC_{\tBunP}^{\leq \mu}\simeq 
\hl_{M/[M,M]!}\fr'_!(\CA^\mu_M\tboxtimes\IC_{\tBunP}^{\leq 0})^r.
$$

Now the map $\fr'$ factors into the projection of the left hand factor
$$
\H_M\underset{\Bun_M\times C}{\times}\tBunP \to 
\H_{M/[M,M]}\underset{\Bun_{M/[M,M]}\times C}{\times}\tBunP 
$$
followed by the projection of the right hand factor
$$
\H_{M/[M,M]}\underset{\Bun_{M/[M,M]}\times C}{\times}\tBunP 
\stackrel{\fr}{\to}
\H_{M/[M,M]} \underset{\Bun_{M/[M,M]}\times C}{\times}\bBunP.
$$
Thus we have an isomorphism
$$
\fr'_!(\CA^\mu_M\tboxtimes\IC_{\tBunP}^{\leq 0})^r
\simeq
\sum_{\nu\in\cowts_M} 
(\IC_{\bBunP}^{\leq 0}
\oplus\CI^{\leq 0})\otimes
\Hom_\chT(V^\nu_\chT,V^\mu_\chM)[\langle2\chrho_M,\nu\rangle]
$$
where as before
$$
\fr_!\IC_{\tBunP}^{\leq 0}\simeq \IC_{\bBunP}^{\leq 0}
\oplus
\CI^{\leq 0}
$$
where $\CI^{\leq 0}$ is isomorphic to a direct sum
of shifts of sheaves of the form 
$$
\IC_{\bBunP}^{\leq\Theta}, 
\mbox{ for pairs 
$\Theta=(0,\fU(\theta^\pos))$, with 
$\theta^\pos\in\poscowts_{G,P}\setminus\{0\}$.}
$$
Finally, applying the modification $\hl_{M/[M,M]!}$ with twist $r(\mu)$
to the above isomorphism, we obtain
an isomorphism
$$
\fr_!\IC_{\tBunP}^{\leq \mu} 
\simeq
\sum_{\nu\in\cowts_M} 
(\IC_{\bBunP}^{\leq r(\mu)}
\oplus \CI^{\leq \mu})\otimes
\Hom_\chT(V^\nu_\chT,V^\mu_\chM)[\langle2\chrho_M,\nu\rangle].
$$
Here we write $\CI^{\leq \mu}$ for the result of applying the modification $\hl_{M/[M,M]!}$ with twist $r(\mu)$
to $\CI^{\leq 0}$. 
Clearly the modification $\hl_{M/[M,M]!}$ takes strata to strata so we conclude that
$\CI^{\leq \mu}$ is isomorphic to a direct sum
of shifts of sheaves of the form 
$$\IC_{\bBunP}^{\leq\Theta}, 
\mbox{ for pairs 
$\Theta=(\theta,\fU(\theta^\pos))$.}
$$
\end{proof}

Note that the proof actually shows that $\CI^{\leq \mu}$
is isomorphic to a direct sum
of shifts of sheaves of the form 
$$
\IC_{\bBunP}^{\leq\Theta}, 
\mbox{ for pairs 
$\Theta=(0,\fU(\theta^\pos))$, with 
$\theta^\pos\in\poscowts_{G,P}\setminus\{0\}$,}
$$
and so in particular is irrelevant, but we shall have no need for this.

Combining the formulas given by Theorem~\ref{thmtBunPconv} and the
preceding lemma, 
we may write the right hand side of equation~\ref{eqnstep1} in the form
\begin{equation}\label{eqnstep1right}
\fr_! H^\lambda_G (\IC_{\tBunP}^{\leq 0})
\simeq
\sum_{\theta\in\cowts_{M/[M,M]}}
\sum_{
{\mu\in\cowts_M},
{r(\mu)=\theta} 
}
\IC_{\bBunP}^{\leq\theta} \otimes \Hom_{\chT} (V_\chT^\mu,V_\chG^\lambda)
[\langle2\chrho_M,\mu\rangle]\oplus\CJ
\end{equation}
where $\CJ$ 
is isomorphic to a direct sum
of shifts of sheaves of the form 
$$\IC_{\bBunP}^{\leq\Theta}, 
\mbox{ for pairs 
$\Theta=(\theta,\fU(\theta^\pos))$.}
$$

Finally, comparing the left hand side (equation~\ref{eqnstep1left})
and the right hand side (equation~\ref{eqnstep1right}), and noting that
$\IC_{\bBunP}^{\leq\theta}$ is not irrelevant, 
we conclude
that 
$$
H^\lambda_G (\IC_{\bBunP}^{\leq 0})\simeq
\sum_{\theta\in\cowts_{M/[M,M]}}
\sum_{
{\mu\in\cowts_M},
{r(\mu)=\theta} 
}
\IC_{\bBunP}^{\leq\theta} \otimes \Hom_{\chT} (V_\chT^\mu,V_\chG^\lambda)
[\langle2\chrho_M,\mu\rangle]\oplus\CM
$$
where $\CM$ is is isomorphic to a direct sum
of shifts of sheaves of the form 
$$\IC_{\bBunP}^{\leq\Theta}, 
\mbox{ for pairs 
$\Theta=(\theta,\fU(\theta^\pos))$.}
$$

\vspace{1em}
{\em Step 2.}
Now we shall show that $\CM$ is in fact zero.
To do this, we shall show that its restriction to each stratum of $\bBunP$
is zero.

Let 
$\Phi$ be a pair
 $(\phi,\fU(\phi^\pos))$,
with $\phi\in\cowts_{M/[M,M]}$,
and $\phi^\pos\in\poscowts_{G,P}$.
Let $H^\lambda_G (\IC_{\bBunP}^{\leq 0})_\Phi$ 
be the restriction of 
$H^\lambda_G (\IC_{\bBunP}^{\leq 0})$ 
to the stratum $\bBunP^\Phi$.
For $\theta\in\cowts_{M/[M,M]}$, let $\CA_\Phi^\theta$ be the 
restriction of $\IC_{\bBunP}^{\leq\theta}$ to the stratum $\bBunP^\Phi$,
and let $\CM_\Phi$ be the restriction of $\CM$.
Note that by step 1, \cite[Theorem 7.3]{BFGMsel02} and Lemma~\ref{lemcyclic} below, 
all of the restrictions are locally constant.

We shall calculate $H^\lambda_G (\IC_{\bBunP}^{\leq 0})_\Phi$ 
in two different ways and compare the results.

On the one hand, by Step 1, we have
\begin{equation}\label{eqnstep1result}
H^\lambda_G (\IC_{\bBunP}^{\leq 0})_\Phi
\simeq
\sum_{\theta\in\cowts_{M/[M,M]}}
\sum_{
{\mu\in\cowts_M},
{r(\mu)=\theta} 
}
\CA^\theta_\Phi\otimes \Hom_{\chT} (V_\chT^\mu,V_\chG^\lambda)
[\langle2\chrho_M,\mu\rangle]\oplus\CM_\Phi
\end{equation}

On the other hand, let us return to the definition of the convolution,
and consider the diagram
$$
\begin{array}{ccccc}
\bBunP & \stackrel{\hlG}{\leftarrow} & \H_G\underset{\Bun_G\times C}{\times}\bBunP^{\leq 0} 
& \stackrel{{}\hrG}{\to} & \bBunP^{\leq 0}\\
\downarrow & & \downarrow & & \downarrow \\
\Bun_G & \stackrel{\hlG}{\leftarrow} & \H_G & \stackrel{{}\hrG}{\to} & \Bun_G
\end{array}
$$
Recall that by definition
$$
H^\lambda_G (\IC_{\bBunP}^{\leq 0})=
{}{\hlG}_!(\CA^\lambda_G\tboxtimes\IC_{\bBunP}^{\leq 0})^r
$$ 
where $(\CA^\lambda_G\tboxtimes\IC_{\bBunP}^{\leq 0})^r$
is the twisted product defined with respect to ${}\hrG$, and $\CA_G^\lambda$
is the simple spherical sheaf on the fibers of $\hrG$ corresponding to $\lambda$.

To calculate 
$H^\lambda_G (\IC_{\bBunP}^{\leq 0})_\Phi$, consider
the inverse image ${h_G^{\leftarrow}}^{ -1}(\bBunP^\Phi)$.
Projecting 
along $\hrG$,
we may decompose the inverse image into a union of locally closed
substacks
$$
{h_G^{\leftarrow}}^{ -1}(\bBunP^\Phi)\simeq
\bigsqcup_{\xi\in\poscoroots_{G,P}} \CS_{P,\phi-\xi}^{\lambda} \timesBunP
\bBunP^{(\xi,\fU(\phi^\pos))}.
$$
Projecting each piece back along 
$\hlG$, we arrive at a spectral sequence for 
$H^\lambda_G (\IC_{\bBunP}^{\leq 0})_\Phi$
with  
$E_2$ term 
$$
\sum_{\xi\in\poscoroots_{G,P}}
\sum_{
{\mu\in\cowts_M},
{r(\mu)=\phi-\xi} 
}
\CA^0_{(\xi,\fU(\phi^\pos))}\otimes \Hom_{\chT} (V_\chT^\mu,V_\chG^\lambda)
[\langle2\chrho_M,\mu\rangle]
$$
In fact, the spectral sequence degenerates here for reasons of parity,
but we shall not need this. What we do need is the following cyclicity. 

\begin{lem}\label{lemcyclic}
Let $\Psi$ be a pair
$(\psi,\fU(\psi^\pos))$,
with $\psi\in\cowts_{M/[M,M]}$,
and $\psi^\pos\in\poscowts_{G,P}$.
Let $\theta\in\cowts_{M/[M,M]}$.
Then
$
\CA^0_{(\psi,\fU(\psi^\pos))}\simeq\CA^\theta_{(\psi+\theta,\fU(\psi^\pos))}.
$
\end{lem}

\begin{proof}
The modification
$$(c,\CP_G,\CP_{M/[M,M]},\sigma)\mapsto
(c,\CP_G,\CP_{M/[M,M]}(\theta\cdot c),\sigma).
$$ 
defines an isomorphism 
$
\bBunP\risom\bBunP
$
which restricts to an isomorphism
$$
\bBunP^{(\psi,\fU(\psi^\pos))}\risom\bBunP^{(\psi+\theta,\fU(\psi^\pos))}.
$$
\end{proof}

We apply the lemma with $\psi=\xi$, $\psi^\pos=\phi^\pos$,
and make the substitution $\theta=\phi-\xi$, to write the $E_2$ term
\begin{equation}\label{eqnstep2}
\sum_{\phi-\theta\in\poscoroots_{G,P}}
\sum_{
{\mu\in\cowts_M},
{r(\mu)=\theta} 
}
\CA^\theta_{(\phi,\fU(\phi^\pos))}\otimes \Hom_{\chT} (V_\chT^\mu,V_\chG^\lambda)
[\langle2\chrho_M,\mu\rangle]
\end{equation}

Comparing our two calculations (equations~\ref{eqnstep1result} and \ref{eqnstep2}),
we conclude by a dimension count that $\CM_\Phi$ must be zero.
\end{proof}
\end{subsection}


\begin{subsection}{Convolution on $\bZcan$}

\begin{thm}\label{thmbZcanconv}
For $\lambda\in\domcowts_G$, there is an isomorphism
$$
 H^\lambda_G (\IC_{\bZcan}^{\leq 0})\simeq
\sum_{\theta\in\cowts_{M/[M,M]}}
\sum_{
{\mu\in\cowts_M},
{r(\mu)=\theta}
}
\IC_{\bZcan}^{\leq\theta} \otimes \Hom_{\chT} (V_\chT^\mu,V_\chG^\lambda)
[\langle2\chrho_M,\mu\rangle].
$$
\end{thm}

\begin{proof}
By Proposition~\ref{propbZcanCart}, for $\theta\in\cowts_{M/[M,M]}$,
we have
$$
\fk^*\IC_{\bBunP}^{\leq \theta}\simeq
\IC_{\bZcan}^{\leq \theta},
$$
Clearly the pullback $\fk^*$ commutes with convolution
$$
H^\lambda_G (\fk^*\IC_{\bBunP}^{\leq \theta})\simeq
\fk^*H^\lambda_G (\IC_{\bBunP}^{\leq \theta}).
$$
Thus by Theorem~\ref{thmbBunPconv}, we conclude
$$
\begin{array}{lcl}
H^\lambda_G (\IC_{\bZcan}^{\leq 0}) & \simeq & 
H^\lambda_G (\fk^*\IC_{\bBunP}^{\leq 0})  \\
& \simeq & \fk^* H^\lambda_G (\IC_{\bBunP}^{\leq 0}) \\
& \simeq &
\sum_{\theta\in\cowts_{M/[M,M]}}
\sum_{
{\mu\in\cowts_M},
{r(\mu)=\theta}}
\fk^*\IC_{\bBunP}^{\leq\theta} \otimes \Hom_{\chT} (V_\chT^\mu,V_\chG^\lambda)
[\langle2\chrho_M,\mu\rangle] \\
& \simeq &
\sum_{\theta\in\cowts_{M/[M,M]}}
\sum_{
{\mu\in\cowts_M},
{r(\mu)=\theta}}
\IC_{\bZcan}^{\leq\theta} \otimes \Hom_{\chT} (V_\chT^\mu,V_\chT^\lambda)
[\langle2\chrho_M,\mu\rangle].
\end{array}
$$
\end{proof}

\end{subsection}


\begin{subsection}{Convolution on $\Z$}
Recall the map of coweight lattices
$$q:
\cowts_M\stackrel{r}{\to}\cowts_{M/[M,M]}\stackrel{p}{\to}\cowts_{A_0}.
$$

\begin{thm}\label{thmZconv}
For $\lambda\in\domcowts_G$, there is an isomorphism
$$
H^\lambda_G (\IC_{\Z}^{\leq 0})\simeq
\sum_{\kappa\in\cowts_{A_0}}
\sum_{
{\mu\in\cowts_{T}},
{q(\mu)=\kappa}
}
\IC_{\Z}^{\leq\kappa} \otimes \Hom_{\chT} (V_\chT^\mu,V_\chG^\lambda)
[\langle2\chrho_M,\mu\rangle].
$$
\end{thm}

\begin{proof}
By Corollary~\ref{corsmall}, for $\theta\in\cowts_{M/[M,M]}$, we have 
$$
\fp_!\IC_{\bZcan}^{\leq\theta}\simeq\IC_{\Zcan}^{\leq p(\theta)},
$$
By Corollary~\ref{corsembed}, for $\kappa\in\cowts_{A_0}$,
we have 
$$
\fs_!\IC_{\Zcan}^{\leq \kappa}\simeq\IC_{\Z}^{\leq \kappa}.
$$
Clearly the pushforwards $\fp_!$ and $\fs_!$ commute with convolution
$$
 H^\lambda_G (\fs_!\fp_!\IC_{\bZcan^{\leq 0}})\simeq
\fs_!\fp_! H^\lambda_G (\IC_{\bZcan^{\leq 0}}).
$$
Thus by Theorem~\ref{thmbZcanconv}, we conclude
$$
\begin{array}{lcl}
H^\lambda_G (\IC_{\Z}^{\leq 0}) & \simeq & 
H^\lambda_G (\fs_!\fp_!\IC_{\bZcan}^{\leq 0})  \\
& \simeq & \fs_!\fp_! H^\lambda_G (\IC_{\bZcan}^{\leq 0}) \\
& \simeq &
\sum_{\theta\in\cowts_{M/[M,M]}}
\sum_{
{\mu\in\cowts_M},
{r(\mu)=\theta} 
}
\fs_!\fp_!\IC_{\bZcan}^{\leq\theta} \otimes \Hom_{\chT} (V_\chT^\mu,V_\chG^\lambda)
[\langle2\chrho_M,\mu\rangle] \\
& \simeq &
\sum_{\kappa\in\cowts_{A_0}}
\sum_{
{\mu\in\cowts_T},
{q(\mu)=\kappa} 
}
\IC_{\Z}^{\leq \kappa} \otimes \Hom_{\chT} (V_\chT^\mu,V_\chG^\lambda)
[\langle2\chrho_M,\mu\rangle].
\end{array}
$$
\end{proof}

\end{subsection}

\end{section}


\begin{section}{Complements}

For our application~\cite{GNcombo06}, we need a slight modification of our main result.
As usual, let $X$ be an affine horospherical $G$-variety with dense 
$G$-orbit $\openX\subset X$ and generic stabilizer $S\subset G$.
Let $S^0$ be the identity component of $S$, and let $\pi_0(S)$ be the component
group $S/S^0$.

For a scheme $\CS$, we write $C_\CS$ for the product $\CS\times C$.
For an $\CS$-point $(c,\CP_G,\sigma)$ of the ind-stack $Z$, 
the section $\sigma$ defines a reduction
of the $G$-bundle $\CP_G$ to an $S$-bundle $\CP'_S$ over an open
subscheme $C_\CS'\subset C_\CS$ which is the complement $C_\CS\setminus\D$
of a subscheme $\D\subset C_\CS$ which is finite and flat over $\CS$. 
By induction, the $S$-bundle $\CP'_S$ defines a
$\pi_0(S)$-bundle 
over $C_\CS'$. 
We call this the generic $\pi_0(S)$-bundle
associated to the point $(c,\CP_G,\sigma)$.

We define $'Z\subset Z$ to be the ind-substack whose $\CS$-points
$(c,\CP_G,\sigma)$ have the property that for every geometric point $s\in\CS$,
the restriction of the associated generic $\pi_0(S)$-bundle to $\{s\}\times C\subset C_\CS$
is trivial. It is not difficult (see \cite{GNcombo06}) to show that $'Z$ is closed in $Z$.
Observe that we have a short exact sequence
$$
0\to\cowts_{A_0}\to\cowts_A\to S/S^0\to 0.
$$
Thus for $\kappa\in\cowts_{A_0}$, it makes sense to consider
the locally closed substack $'Z^{\kappa}\subset {}'Z$
and its closure $'Z^{\leq\kappa}\subset {}'Z$.
Observe as well that from the fibration $S\to G\to G/S$, we have an exact sequence
$$
\pi_1(G)\to\pi_1(\openX)\to \pi_0(S).
$$
Thus for $\lambda\in\domcowts_G$, 
we have the convolution functor 
$$
H^\lambda_G :\Sh('Z)\to\Sh('Z).
$$ 

The same arguments show that our main result holds equally well in this context.

\begin{thm}\label{thmZconv}
For $\lambda\in\domcowts_G$, there is an isomorphism
$$
H^\lambda_G (\IC_{'\Z}^{\leq 0})\simeq
\sum_{\kappa\in\cowts_{A_0}}
\sum_{
{\mu\in\cowts_{T}},
{q(\mu)=\kappa}
}
\IC_{'\Z}^{\leq\kappa} \otimes \Hom_{\chT} (V_\chT^\mu,V_\chG^\lambda)
[\langle2\chrho_M,\mu\rangle].
$$
\end{thm}

\end{section}


\bibliographystyle{alpha}
\bibliography{ref}

\end{document}